\documentclass{amsart}
\usepackage[parfill]{parskip} 
\usepackage{amsmath,amssymb,amsthm}
\usepackage[pdftex]{hyperref}
\usepackage[margin=2.5cm]{geometry}
\usepackage{euscript}
\usepackage{tikz-cd}
\usepackage{appendix}
\usepackage{comment}
\usepackage{mathtools}
\usepackage{graphicx, pgf} 
\usepackage[scr=boondox]{mathalpha} 

\hypersetup{
  colorlinks,
  citecolor=blue,
  filecolor=black,
  linkcolor=blue,
  urlcolor=black
}

\numberwithin{equation}{section}

\theoremstyle{plain}
\newtheorem{theorem}{Theorem}[section]
\newtheorem*{theorem*}{Theorem}
\newtheorem{lemma}[theorem]{Lemma}
\newtheorem{proposition}[theorem]{Proposition}
\newtheorem{corollary}[theorem]{Corollary}
\newtheorem{example}[theorem]{Example}
\newtheorem{innercustomthm}{Theorem}
\newenvironment{customthm}[1]
  {\renewcommand\theinnercustomthm{#1}\innercustomthm}
  {\endinnercustomthm}

\theoremstyle{definition}
\newtheorem{definition}[theorem]{Definition}

\newtheorem{remark}[theorem]{Remark}

\newcommand{\updownarrows}{\mathbin\uparrow\hspace{-.1em}\downarrow}
\newcommand{\downuparrows}{\mathbin\downarrow\hspace{-.1em}\uparrow}
\newcommand{\e}{\mathsf e} 
\newcommand*{\parallelogramm}[1][]{%
  \pgfpicture\pgfsetroundjoin
    \pgftransformxslant{.3}%
    \pgfpathrectangle{\pgfpointorigin}{\pgfpoint{.35em}{.35em}}%
    \pgfusepath{stroke,#1}%
  \endpgfpicture}
\def\epsilon{\varepsilon}
\def\phi{\varphi}
\def\uptriangle{\vartriangle}
\def\emptyset{\varnothing}
\def\IMP{\Rightarrow}
\def\IFF{\Leftrightarrow}
\def\downtriangle{\triangledown}
\DeclareMathOperator{\End}{End}
\DeclareMathOperator{\Mod}{Mod}
\DeclareMathOperator{\Clo}{Clo}
\DeclareMathOperator{\Th}{Th}

\DeclareMathOperator{\Str}{Str}

\DeclareMathOperator{\set}{set}
\DeclareMathOperator{\fin}{fin}

\DeclareMathOperator{\Max}{Max}

\begin{document}
\title[Birkhoff-style Theorems Through Infinitary Clone Algebras]{Birkhoff-style Theorems Through Infinitary Clone Algebras}
\author[A. Bucciarelli]{Antonio Bucciarelli}
\author[P.-L. Curien]{Pierre-Louis Curien}
\author[A. De Faveri]{Arturo De Faveri}
\author[A. Salibra]{Antonino Salibra}

\address{IRIF\\
CNRS and Universit\'e Paris Cit\'e\\
8 Place Aur\'elie Nemours, 75013 Paris, France
}

\begin{abstract}
  Building upon the classical article ``Representing varieties of algebras by algebras'' by W. D. Neumann, 
  we revisit the famous Birkhoff's HSP theorem in the light of infinitary algebra.  
\end{abstract}

\maketitle

\section{Introduction}
Just as monoids serve to represent the composition of functions from a set $A$ to itself, 
clones aim at modelling the composition of functions $A^n \to A$ of arbitrary finite arity $n \in \omega$.
Clones were originally defined as sets of functions on $A$ that contain all the projections and are closed under composition, and can be abstracted as many-sorted algebras \cite[Definition 10.2]{ALV3}
(for a historical account of clones see \cite[p. 126]{ALV3}). 
Clones hold significant importance in universal algebra, enabling the study of varieties in a way independent of their presentation (i.e.\ of their similarity type and equations).
However, the study of clones presents challenges, primarily due to their $\omega$-sorted structure, featuring a distinct sort for each finite arity.

One way to circumvent this problem is to model the composition of functions $A^\omega \to A$ of fixed arity $\omega$
and to recover the case of finitary operations as a particular case: any $n$-ary operation, $n \in \omega$, can be viewed as an $\omega$-ary operation depending only on the first $n$ variables.  
The idea of developing a theory of clones of infinitary operations to encode clones as one-sorted algebras was exploited in \cite{BS22}, where the authors introduced a class of algebras (labeled \emph{clone algebras}) 
with a countable number of nullary operations $\e_0, \e_1, \e_2, \ldots$ and a $(n+1)$-ary operation $q_n$ for every $n \in \omega$. 
Roughly speaking, $q_n(f, g_1,\ldots, g_n)$ represents the composition of $f$ with $g_1,\ldots,g_n$ in the first $n$ coordinates of $f$. 
Another approach, which comes at the expense of departing from the realm of classical (finitary) universal algebra, 
originated in the work of Walter D. Neumann. 
In \cite{neu70} he introduced \emph{$\aleph_0$-clones}, a class of infinitary algebras
whose type 
consists, besides nullary operations $\e_0, \e_1, \e_2, \ldots$, of a single $\omega$-ary operation $q$ representing infinitary composition. 
Neumann also proved a Cayley-style representation theorem: 
just as any monoid is a transformation monoid of its underlying set via the action given by left multiplication, 
every $\aleph_0$-clone is (isomorphic to) a clone of $\omega$-ary functions on its underlying set. 



In this paper, we augment the type of $\aleph_0$-clones with a nullary operation  $f$ for each $f \in \tau$, for a given a similarity type $\tau$ of infinitary operation symbols.
We define an equational class of infinitary algebras, called \emph{infinitary clone $\tau$-algebras}, over this augmented type. 
In infinitary clone $\tau$-algebras the nullary operations  $\e_0, \e_1, \e_2, \ldots$ abstract the role played by variables in free algebras and projections in clones, and the operation $q$ gives a unified way to handle term-for-variable substitution and functional composition. 
Among {infinitary} clone $\tau$-algebras, a special role is played by the \emph{functional} ones. 
The universe of a functional infinitary clone $\tau$-algebra $\mathcal{C}$ is made of $\omega$-ary functions that include the operations of a given type $\tau$-algebra $\mathbf{A}$,  
called the \emph{value domain} of $\mathcal{C}$. 

Given a type $\tau$, we consider the interaction between two classes of algebras. 
At the bottom level, we have the $\tau$-algebras ordinarily studied; 
at the top level, there are the {infinitary} clone $\tau$-algebras. 
The goal is to study the interplay between these two levels and to apply techniques from the higher level to derive results in the lower one. 
More concretely, this interplay develops as follows. 

If $\mathsf{K}$ is a class of $\tau$-algebras, we define $\mathsf{K}^{\uptriangle}$ as the class of the functional infinitary clone $\tau$-algebras (up to isomorphism) with value domains in $\mathsf{K}$; 
in the converse direction, if $\mathsf{H}$ is a class of infinitary clone $\tau$-algebras, 
we define $\mathsf{H}^{\downtriangle}$ as the class of $\tau$-algebras (up to isomorphism) which appear as value domains of some elements of $\mathsf{H}$. 
We prove that in both cases, starting from a variety in the lower (resp. higher) level, the resulting class is again a variety in the higher (resp. lower) one. 

The machinery described so far allows for an enrichment of a pivotal result in universal algebra. 
Birkhoff characterised those classes of algebras that can be defined by a set of identities as those closed under products, subalgebras, and homomorphic images (called varieties). 
We extend this characterisation, in the infinitary case, with a condition expressed in the language of clone algebras.
\begin{customthm}{6.7}
  \label{fakebirk}
For a class $\mathsf{K}$ of $\tau$-algebras, the following are equivalent:
\begin{enumerate}
\item $\mathsf{K}$ is a variety of $\tau$-algebras; 
\item $\mathsf{K}$ is an equational class;
\item $\mathsf{K}=\mathsf{K}^{\uptriangle \downtriangle}$ and $\mathsf{K}^{\uptriangle}$ is a variety of infinitary clone $\tau$-algebras.
\end{enumerate}
\end{customthm}
The proof of the known equivalence $(1) \IFF (2)$ is new.  
We also give examples showing that the two conditions in $(3)$ are both necessary.

As a corollary, we obtain the following result concerning finitary algebras. 
Given a finitary type $\rho$, let $\rho_\bullet$ be the type obtained by considering every operation symbol in $\rho$ as $\omega$-ary.
Accordingly, given a class $\mathsf{V}$ of $\rho$-algebras, let $\mathsf{V}^{\top}$
be the class of $\rho_\bullet$-algebras obtained by considering the operations of $\mathsf{V}$ formally infinitary.
\begin{customthm}{7.8}
  Let $\rho$ be a finitary type
  and $\mathsf{V}$ be a class of $\rho$-algebras. 
  The following are equivalent:
  \begin{enumerate}
  \item $\mathsf{V}$ is a variety of $\rho$-algebras;
  \item $\mathsf{V}$ is an equational class;
  \item $\mathsf{V}^{\top}$ is a variety of $\rho_\bullet$-algebras; 
  \item $\mathsf{V}^{\top}=\mathsf{(V^{\top})}^{\uptriangle \downtriangle}$ and $\mathsf{(\mathsf{V}^{\top})}^{\uptriangle}$ is a variety of infinitary clone $\rho_\bullet$-algebras.
  \end{enumerate}
  \end{customthm}
We apply the same technique to the study of topological
variants of Birkhoff Theorem, as initiated by Bodirsky-Pinsker \cite{BP15}, and generalised recently by Schneider \cite{Sc17} and Gehrke-Pinsker \cite{GP18}. 
These authors provide a Birkhoff-type characterisation of all those members of the pseudovariety generated by a given algebra. 
By endowing infinitary clone algebras with a suitable topology, we obtain a simple topological version of Birkhoff Theorem for infinitary algebras, which turns out to subsume the results cited above. 

As a final note in this introduction, we address the question why we have chosen to employ infinitary clone algebras rather than clone algebras in the present work 
(indeed, the intended models of both structures are clones of infinitary operations).
As a matter of fact, clone algebras fail to completely grasp the structure of clones of infinitary operations:  
every operator $q_n$ can only compose a finite number of infinitary operations.  
As a consequence, using clone algebras instead of infinitary clone algebras at the upper level of the interplay described above 
deeply affects the very nature of the algebraic structures at the lower level. 
This issue can be resolved by using an \emph{ad hoc} definition for these algebraic structures, as 
initially explored by the fourth author in a preliminary version of this work \cite{S22}.
Our approach in this work avoids all the unnecessary complications of that preliminary version.  

We now outline the plan of this work. 
In Section \ref{sec:prelim}, we present some preliminary notions. 
Section \ref{sec:neu} introduces the main object of study: infinitary clone algebras. 
Section \ref{sec:cats} demonstrates how terms, equations, and free algebras, in the infinitary case,
can be encoded within the language of infinitary clone algebras.
The technique of associating a variety of infinitary algebras with a variety of infinitary clone algebras, and vice versa, is described in Section \ref{sec:up-down}. 
Section \ref{sec:gb} contains the core of this work: Theorem \ref{fakebirk}. 
Section \ref{sec:tafromalg} substantiates the claim that these results can be applied to finitary algebras, thus giving an enhanced version of Birkhoff Theorem in that case too. 
Section \ref{sec:top} of the paper is dedicated to the topological version of Birkhoff Theorem in our setting. 
The final section suggests some future directions of work.  

\section{Preliminaries}
\label{sec:prelim}
The notation and terminology of this paper are pretty standard. 
For concepts, notations and results not covered hereafter, the reader is
referred to \cite{BS81,mac87} for universal algebra, to \cite{SZ86, T93, L06} \cite[Chapter 10]{ALV3} for the theory of clones.

If $s \in A^\omega$, we write $\set(s)$ for $\{s_i : i \in \omega\}$.  
We denote by $A^\omega_{\fin}$ the set $\{s \in A^\omega : |\set(s)| < \omega\}$. 
We shall sometimes use a bold symbol $\boldsymbol{x}$ for a sequence $(x_0,x_1, \ldots)$ of length $\omega$. 
If $f:A\to B$ is a function, we denote by $f^\omega$ the function $A^\omega \to B^\omega$ such that $f^\omega(s)= (f(s_i) : i\in\omega)$ for every $s\in A^\omega$.

By a \emph{finitary operation} on a set $A$ we  mean  a function $f:A^n\to A$ for some natural number $ n \in \omega$. 
We denote by $O_A$ the set of finitary operations on a set $A$. 
If $F\subseteq O_A $, then $F^{(n)}:= \{f:A^n\to A \, | \,  f\in  F\}$. 

By an $\omega$-\emph{ary operation} on $A$ we mean a function $f: A^\omega \to A$. 
The set of $\omega$-ary operations on $A$ is denoted by $O_A^{(\omega)}$.
We shall often call an operation of arity $\omega$ just an ``infinitary operation''.
 

A type $\sigma$ is a set of operation symbols, coming each with an arity $ \le \omega$. 
We shall say that $\sigma$ is
\begin{itemize}
  \item \emph{finitary} if each symbol has arity $< \omega$; 
  \item \emph{infinitary} if there is at least one symbol of arity $\omega$; 
  \item \emph{homogeneous} if every function symbol has the same arity $\omega$. 
\end{itemize}
For every $\kappa \le \omega$, we denote by $\sigma_\kappa$ the set of function symbols of arity $\kappa$. 

A $\sigma$-algebra is a tuple $\mathbf A=(A,f^\mathbf A)_{f\in\sigma}$ such that $f^\mathbf A:A^\kappa \to A$ is an operation of arity $\kappa \leq \omega$ for every $f\in\sigma_\kappa$. 
We denote by $\mathsf{Alg}_{\sigma}$ the class of all $\sigma$-algebras. 
Algebras of type $\sigma$ and homomorphisms between them form a category. 
An infinitary algebra is an algebra over an infinitary type. 
We refer the reader to \cite{slo} for an account of infinitary algebras.

A $\sigma$-algebra $\mathbf{A}$ is \emph{minimal} if it does not have any proper subalgebra. 
Every $\sigma$-algebra $\mathbf{A}$ admits a unique minimal subalgebra, 
which is the subalgebra of $\mathbf{A}$ generated by the interpretation of the $0$-ary symbols of $\sigma$ in $\mathbf{A}$. 
In particular the minimal subalgebra of $\mathbf{A}$ is empty iff $\sigma$ has no $0$-ary symbol. 
The homomorphic image of any minimal algebra is minimal. 

Let $X$ be a set of variables. 
\begin{definition}
  \label{def_terms}
  We define sets $A_{\mu}$ by transfinite induction on $\mu \in \mathsf{Ord}$.
  \begin{enumerate}
    \item $A_{0} = X$; 
    \item $A_{\mu} = A_{\nu} \cup \{f(t_0, \ldots, t_{n-1}) : t_i \in A_\nu, f \in \sigma_n, n \in \omega\} \cup  \{f(t_0, t_1, \ldots) : t_i \in A_{\nu}, f \in \sigma_\omega\}$ for $\mu = \nu +1$;
   \item $A_{\mu} = \bigcup_{\nu < \mu} A_{\nu} $ for $\mu$ limit.
\end{enumerate}
We define $T_{\sigma}(X):=A_{\omega_1}$ as the set of $\sigma$-\emph{terms} over $X$. 
If $X=\emptyset$ we simply write $T_\sigma$; if $\sigma_0 = \emptyset$, then $T_\sigma = \emptyset$.  
\end{definition}
Notice that this inductive definition yields a well-order which is the usual subterm ordering on $T_{\sigma}(X)$, that is the smallest containing
the order $t_i < f(t_0, \ldots, t_n, \ldots)$.

The initial object $\mathbf{I}$ in the category $\mathsf{Alg}_{\sigma}$, whose elements are the closed terms over $\sigma$, is non-empty iff there are symbols of arity $0$ in the signature. 
The initial algebra is minimal since it is generated by the $0$-ary symbols.

If $\mathbf{A}$ is a $\sigma$-algebra we write $\Th(\mathbf{A})$ for the set of identities, i.e. pairs of $\sigma$-terms, satisfied by $\mathbf{A}$. 
If $\mathsf{K}$ is a class of $\sigma$-algebras, we write $\Th(\mathsf{K})$ for $\bigcap\{\Th(\mathbf{A}) : \mathbf{A} \in \mathsf{K}\}$. 
The set $\Th(\mathbf{A})$ (resp. $\Th(\mathsf{K}$)) is called the \emph{equational theory} of $\mathbf{A}$ (of $\mathsf{K}$ resp.). 
If $\Sigma$ is a set of identities, we write $\Mod(\Sigma)$ for the class of models of $\Sigma$, i.e those $\sigma$-algebras satisfying each identity in $\Sigma$. 
A class of similar algebras $\mathsf{K}$ such that $\mathsf{K}=\Mod(\Sigma)$ for some $\Sigma$ is called \emph{equational class}. 
If $\mathsf{K}$ is a class of $\sigma$-algebras, then we denote by $\mathbf F_\mathsf{K}(X)$, if it exists, the free $\mathsf{K}$-algebra over $X$.  
If $\mathsf{K}$ is an equational class, $\mathbf F_\mathsf{K}(X)$ always exists in $\mathsf{K}$.

Closure of a class of similar algebras under homomorphic images, (finite, countable) products, subalgebras and isomorphic images is denoted by $\mathbb{H}$, $\mathbb{P}$ ($\mathbb{P}_{\fin}$, $\mathbb{P}_{\omega}$), $\mathbb{S}$ and $\mathbb{I}$, respectively. 
A class of similar algebras $\mathsf{K}$ such that $\mathsf{K}=\mathbb{HSP}(\mathsf{K})$ is called \emph{variety};
$\mathbb{V}(\mathsf{K})$ is the smallest variety containing $\mathsf{K}$.


Infinitary operations appear naturally in various contexts. 
Every semilattice, lattice, or Boolean algebra $\mathbf{A}$ with countable joins can be seen as algebra with some finitary operations
and an infinitary operation $f : A^\omega \to A$, $f(s) =\sup\{s_i : i \in \omega\}$.
As another example, the opposite category of the category of compact Hausdorff topological spaces is equivalent to an equational class of infinitary algebras, called $\delta$-$\emph{algebras}$. 
The type of $\delta$-algebras contains a finite number of finitary operations and a single $\omega$-ary operation $\delta$ \cite{IS82}.
The class of $\delta$-algebras can be axiomatised by a finite number of identities \cite{MR17}. 

A \emph{clone on a set $S$} is a set of finitary operations containing all the projections and closed under composition. 
A \emph{clone on a $\rho$-algebra $\mathbf S$}, where $\rho$ is a finitary type, is a clone on $S$ containing the operations $f^\mathbf S$ of $\mathbf S$, for $f\in\rho$.
The smallest clone on $\mathbf{S}$ is denoted by $\Clo(\mathbf{S})$. 

\section{Infinitary clone algebras}\label{sec:neu} 
Throughout the paper we reserve the symbol $\tau$ to mean a homogeneous type of infinitary operation symbols such that $\tau \cap \{q, \e_i\}_{i \in \omega}=\varnothing$. 
We denote by $\bar{\tau}$ the type $(q, \e_i, f)_{i \in \omega, f \in \tau}$ where $\e_i,f$ are $0$-ary operation symbols and $q$ is an operation symbol of arity $\omega$. 
Note that the arity of symbols of $\tau$ changed from $\omega$ to $0$. 

The following definition aims at capturing the algebraic structure of the set of operations of an algebra (see Definition \ref{def_func}), the set of terms of a type, and the set of term operations of an algebra (Section \ref{sec:cats}).

\begin{definition}
  \label{def:clo}
  An \emph{infinitary clone $\tau$-algebra} is a $\bar{\tau}$-algebra
  \begin{equation*}
    \mathcal{C}=(C, q^\mathcal{C},\e_i^\mathcal{C},f^\mathcal{C})_{i \in \omega,f\in\tau}\text{,}
  \end{equation*}
   where $\e_i^\mathcal C,f^\mathcal C\in C$ 
   and $q^\mathcal C$ is an $\omega$-ary operation satisfying the identities
   \begin{enumerate}
    \item[(N1)] $q(\e_i,x_0,\dots,x_n,\dots) = x_i;$
    \item[(N2)] $q(x,\e_0,\dots,\e_n,\dots)  = x;$
    \item[(N3)] $q(q(x,\boldsymbol{y}),\boldsymbol{z}) = q(x,q(y_0,\boldsymbol{z}),\dots,q(y_n,\boldsymbol{z}),\dots)$.
   \end{enumerate}
  We denote by $\mathsf{CA}_\tau$ the class of infinitary clone $\tau$-algebras.
  If $\tau=\emptyset$, $\mathcal{C}$ is called \emph{pure}. 
\end{definition}

Pure infinitary clone algebras 
were introduced by Neumann in \cite{neu70} with the name ``(abstract) $\aleph_0$-clones''.
Neumann proved that every $\aleph_0$-clone can be represented by a clone of operations of arity $\omega$ on some set $A$.
This is the content of Proposition \ref{lem_repr}.

If $\mathcal{C} \in \mathsf{CA}_\tau$ and $a\in C$, then $\widetilde{a}^{\mathcal{C}}:C^\omega\to C$ is the infinitary operation defined as follows: 
\begin{equation*}
  \widetilde{a}^{\mathcal C}(s):=q^\mathcal{C}(a,s) \text{,}
\end{equation*}
for every $s\in C^\omega$. 
If $a= f^\mathcal{C}$ for some $f\in \tau$, we write $\widetilde{f^\mathcal{C}}$ for $\widetilde{f^\mathcal{C}}^\mathcal{C}$.

\begin{definition}
  Given $\mathcal{C} \in \mathsf{CA}_\tau$, we define
  the $\tau$-algebra $\mathcal{C}^\downarrow = (C, \widetilde{f^\mathcal{C}})_{f \in \tau}$.
\end{definition}

\begin{lemma}
  \label{rem_homisalso}
  The map $(-)^{\downarrow} : \mathsf{CA}_{\tau} \to \mathsf{Alg}_{\tau}$, $\mathcal{C} \mapsto \mathcal{C}^{\downarrow}$, extends to a functor.
  This functor is the identity on arrows and it preserves isomorphisms. 
  \begin{proof}
    Let $\alpha: \mathcal{C} \to \mathcal{D}$ be a homomorphism. 
    Then we have  $\alpha(\widetilde{f^{\mathcal{C}}}(s)) = \alpha(q^{\mathcal{C}}(f^{\mathcal{C}}, s))= q^{\mathcal{D}}(\alpha(f^{\mathcal{C}}), \alpha^{\omega}(s))=\widetilde{f^{\mathcal{D}}}(\alpha^{\omega}(s))$. 
  \end{proof}
\end{lemma}

In the class of infinitary clone $\tau$-algebras, a special role is played by those that arise as sets of operations of some $\tau$-algebra $\mathbf{A}$. 
 
\begin{definition} 
  \label{def_func}
  Let $\mathbf{A}=(A,f^\mathbf{A})_{f\in\tau}$ be a $\tau$-algebra.
  The \emph{full functional infinitary clone $\tau$-algebra on value domain $\mathbf A$} is the $\bar{\tau}$-algebra 
  \begin{equation*}
    \mathcal{O}_\mathbf{A}^{(\omega)}=(O_A^{(\omega)},\e_i^A,q^A, f^\mathbf A)_{i \in \omega,f\in\tau}
  \end{equation*}
  where $\e_i^A, q^A$ are defined as follows: $\e_i^A(s)=s_i$ and 
  $q^A(g_0,g_1,\dots,g_n,\dots)(s)=g_0(g_1(s),\dots,g_n(s),\dots)\text{,}$
  for all $s \in A^\omega$ and $g_i \in O^{(\omega)}_A$.
  A \emph{functional infinitary clone $\tau$-algebra  on value domain $\mathbf A$}  is a subalgebra of $\mathcal{O}_\mathbf A^{(\omega)}$.
  We denote by $\mathsf{FCA}(\mathbf{A})$ the class of functional infinitary clone $\tau$-algebra on value domain $\mathbf{A}$.
\end{definition}

Note that the interpretation $f^{\mathbf{A}}$ of the $\omega$-ary symbol $f \in \tau$ in $\mathbf{A}$ becomes the interpretation of the nullary symbol $f \in \bar \tau$ in $\mathcal{O}^{(\omega)}_{\mathbf{A}}$. 

The map $\phi: \mathcal{C} \to \mathcal{O}^{(\omega)}_{\mathcal{C}^\downarrow}$, defined by $\phi(a) = \widetilde a^\mathcal{C}$, defines an embedding
(directly adapted from Neumann \cite{neu70}).

\begin{proposition} 
  \label{lem_repr} 
  Every $\mathcal{C} \in \mathsf{CA}_\tau$ is isomorphic to a functional infinitary clone $\tau$-algebra on value domain $\mathcal{C}^\downarrow$.
\end{proposition}

\begin{example} \label{exa:proj} 
  The algebra $\mathcal{P}=(\omega, q^\mathcal{P},\e_i^\mathcal{P})_{i \in \omega}$, where  $\e_i^\mathcal{P}=i$ and $q^\mathcal{P}(i,s)=s_i$ for every $s\in\omega^\omega$,
  is the initial pure infinitary clone algebra.
\end{example}

\begin{example}\label{exa:free} 
  Let $\mathsf{K}$ be an equational class of $\tau$-algebras and $\mathbf F_\mathsf K(X)$ ($\mathbf F$, for short)
  be the free $\mathsf{K}$-algebra over the countable set $X=\{x_0,x_1,\dots\}$ of generators. 
  The $\bar{\tau}$-algebra $\mathcal F_\mathsf K=(F, q^\mathcal{F}, \e_i^{\mathcal F}, f^\mathcal F)_{f\in\tau}$ is an infinitary clone $\tau$-algebra, with operations defined as follows. 
  For every $i \in \omega$ let $\e_i^\mathcal{F}: = x_i \in F$. 
  For every $f \in \tau$, $f^\mathcal F := f^{\mathbf{F}}(\boldsymbol{x}) \in F$.
  Finally, $q^{\mathcal{F}}(a, \boldsymbol{b}) := \phi(a)$ where $\phi$ is the unique endomorphism
  of $\mathbf{F}$ which sends the generators $x_i$ to $b_i$ for every $i \in \omega$. 
\end{example}

The following proposition will be useful in Sections \ref{sec:up-down} and \ref{sec:gb}. 

\begin{proposition}
  \label{prop_rel}
  Let $\mathbf{A}$ be a $\tau$-algebra.
  \begin{enumerate}
  \item \label{lem_val} For every $s \in A^{\omega}$, the function $\bar{s} : (\mathcal{O}^{(\omega)}_{\mathbf{A}})^{\downarrow} \to \mathbf{A}$ mapping $\phi$ to $\phi(s)$ is a homomorphism of $\tau$-algebras. 
  \item  \label{cor_subalg_of_prod}
  Moreover, $(\mathcal{O}^{(\omega)}_{\mathbf{A}})^{\downarrow}$ can be embedded into $\mathbf{A}^{A^\omega}$. 
  \end{enumerate}
  \begin{proof}
    (1) Let $s \in A^\omega$ and $\phi_i \in \mathcal{O}^{(\omega)}_{\mathbf{A}}$. 
    We have
    \begin{align*}
      \bar{s}(\widetilde{f^{\mathbf{A}}}(\phi_0, \ldots, \phi_i, \ldots)) & =\widetilde{f^{\mathbf{A}}}(\phi_0, \ldots, \phi_i, \ldots)(s) \\
      & =q^A(f^{\mathbf{A}}, \phi_0, \ldots, \phi_i, \ldots)(s) \\
      & =f^{\mathbf{A}}(\bar{s}(\phi_0), \ldots, \bar{s}(\phi_i),\ldots)\text{.}
    \end{align*}
    As for (2), the function $\phi \mapsto (\phi(s): s  \in A^\omega)$ is the required embedding. 
  \end{proof}
\end{proposition}

\subsection{Independence and dimension}
We define the notions of independence and dimension in the infinitary clone algebras, abstracting the notion of arity of the finitary operations. 

\begin{definition} 
  \label{def:dim}
  Let $\mathcal C$ be an infinitary clone $\tau$-algebra  and $a \in C$.  
  We say that
  \begin{enumerate}
    \item[(i)] $a$ has \emph{dimension} at most $n$ ($\dim(a) \leq n$) if, for all $s,u\in C^\omega$,  $u_i=s_i$ for every $i=0,\dots,n-1$ implies $\widetilde a^\mathcal C(s) = \widetilde a^\mathcal C(u)$;
    \item[(ii)] $a$ has \emph{dimension $n$} ($\dim(a) = n$) if $n$ is the minimal natural number such that $a$  has dimension at most $n$;
    \item[(iii)] $a$ is \emph{infinite-dimensional} if there is no $n$ such that $\dim(a) \leq n$.
  \end{enumerate}
\end{definition} 

We denote by 
$\mathrm{Fd}(\mathcal C)$ the set of its finite-dimensional elements. 
We say that $\mathcal C$ is \emph{finite-dimensional} if $ C=\mathrm{Fd}(\mathcal C)$.
If every $f \in \tau$ is interpreted in a finite-dimensional element of $\mathcal{C}$,
then $\mathrm{Fd}(\mathcal C)$ is a subalgebra of $\mathcal C$.
The clone algebra introduced in Example \ref{exa:proj} is finite-dimensional.

\begin{lemma} \cite[Section 3]{neu70} 
  \label{lem_neu}
  Let $n\geq 1$,
  $\mathcal{C}$ be an infinitary clone $\tau$-algebra and $a \in C$.  
  The following properties hold:
  \begin{enumerate}
    \item $\dim(a) \le n$ iff
    $q(a,\e_0,\dots,\e_{n-1},\e_0,\e_0,\dots)=a$.
    \item $\dim(a) =0$  iff
    $q(a,\e_0,\e_2,\e_4,\e_6,\dots)=q(a,\e_1,\e_3,\e_5,\e_7,\dots)$.
  \end{enumerate}
\end{lemma}

\subsection{Similarity}
We introduce the top extension of a finitary operation. 
This notion allows us, without loss of generality, to regard every type as a homogeneous type and will be essential in Section \ref{sec:tafromalg}. 
\begin{definition}\label{def:topext}
 Let $f:A^n\to A$ be a finitary operation. The  \emph{top extension $f^\top:A^\omega\to A$} of $f$ is defined as follows:
 $$f^\top(s)=f(s_0,\dots,s_{n-1}),\ \text{for every $s\in A^\omega$}.$$
 We say that $f, g \in O_A$ are \emph{similar}, and we write $f\approx g$, if $f^\top=g^\top$.
\end{definition}

The relation $\approx$ is an equivalence relation. 
If $B$ is an equivalence class of $\approx$, then $B \cap  O^{(n)}_A$ is either empty or a singleton.  
If $B \cap  O^{(n)}_A\neq\emptyset$, then $B \cap  O^{(k)}_A\neq\emptyset$ for every $k\geq n$. 

\begin{lemma}  \cite{BS22}
  Every clone on $A$ is a union of equivalence classes of $\approx$. 
\end{lemma}

The set of all top extensions is denoted by $O_A^\top=\{f^\top : f\in O_A\}$.
If $C$ is a clone on $A$, then $C^\top = \{f^\top : f\in C\}$ is the universe of a finite-dimensional pure functional infinitary clone algebra. 

\begin{lemma}\label{lem:dim} \cite{BS22} Let $\mathcal C\in \mathsf{CA}_{\tau}$ and let $a \in C$. 
  \begin{enumerate}
    \item If $\dim(a)=n$, then there is an $f \in O_C^{(n)}$ such that $\widetilde{a}^{\mathcal{C}}=f^{\top}$. 
    \item $\mathcal C$ is finite-dimensional iff $\mathcal C \simeq D^\top$ for some clone $D$ on $C$.
\end{enumerate}
\end{lemma}
\section{Free infinitary clone algebras}\label{sec:cats}
The following definition is completely analogous to Definition \ref{def_terms}, but is included here for clarity's sake. 
Let $\mathscr{X}$ be a set such that $\mathscr{X} \cap \tau=\emptyset$.

\begin{definition}
  \label{def:meta}
  We define sets $B_{\mu}$ by transfinite induction on $\mu \in \mathsf{Ord}$.
  \begin{enumerate}
    \item $B_{0} = \{\e_0, \e_1, \ldots\}$; 
    \item $B_{\mu} = B_{\nu} \cup \{w(t_0, t_1, \ldots) : t_i \in B_{\nu}, w \in \tau \cup \mathscr{X}\}$ for $\mu = \nu +1$;
   \item $B_{\mu} = \bigcup_{\nu < \mu} B_{\nu} $ for $\mu$ limit.
\end{enumerate}
We define $N_{\bar \tau}(\mathscr{X}):=B_{\omega_1}$ as the set of $\tau$-\emph{metaterms}.  
\end{definition}
 
For $\mathscr{x} \in \mathscr{X}$, let $\dot{\mathscr{x}}=\mathscr{x}(\e_0, \e_1, \ldots)$ and let $\dot{\mathscr{X}}:=\{\dot{\mathscr{x}} : \mathscr{x} \in \mathscr{X}\} $.

\begin{definition}
  \label{def:en}
We define the $\bar{\tau}$-algebra $\mathcal{N}_{\bar \tau}(\mathscr{X}):=(N_{\bar \tau}(\mathscr{X}),q^\mathcal{N},\e_i^\mathcal{N},f^\mathcal{N})_{f\in\tau}$, where
$\e_i^\mathcal{N}=\e_i$;  $f^\mathcal{N}=f(\e_0, \e_1, \ldots)$ for every $f\in\tau$; and $q^\mathcal{N}$ is defined by induction on the first argument as follows:
  \begin{enumerate}
    \item[(i)] $q^\mathcal{N}(\e_i,t_0,\dots,t_k,\dots)=t_i$;
    \item[(ii)] $q^\mathcal{N}(w(t_0,\dots,t_k,\dots),\boldsymbol u) \!= \!
      w(q^\mathcal{N}(t_0,\boldsymbol  u),\dots, q^\mathcal{N}(t_k,\boldsymbol  u),\dots)$ for $w\!\in\!\tau\cup X$. 
  \end{enumerate}
  If $\mathscr{X}=\emptyset$ we write $N_{\bar \tau}$ in place of $N_{\bar \tau}(\emptyset)$. 
\end{definition}


 
We now prove that the set $N_{\bar \tau}(\mathscr{X})$ can be endowed with the structure of infinitary clone $\tau$-algebra and that it is the free infinitary clone $\tau$-algebra.

\begin{proposition}\label{lem:term}
  The $\bar{\tau}$-algebra $\mathcal{N}_{\bar \tau}(\mathscr{X})$
 is the free infinitary clone $\tau$-algebra over the set $\dot{\mathscr{X}}$ of generators in the class $\mathsf{CA}_\tau$.  
 If $\mathscr{X}=\emptyset$, then $\mathcal{N}_{\bar \tau}$
 is initial in $\mathsf{CA}_\tau$.
\begin{proof} 
  We show that $\mathcal{N}_{\bar \tau}(\mathscr{X})$ satisfies the identities of Definition \ref{def:clo}. 
  We leave (N1) and (N2) to the reader.
  As for (N3), consider $q(q(v,\boldsymbol  y),\boldsymbol z)$. 
  We reason
  by induction on $v$. 
  If $v=\e_i$ for some $i$, the proof is trivial. 
  If $v=w(t)$ for $t=(t_0,\dots,t_k,\dots)$ and $w\in \tau\cup \mathscr{X}$:
 \begin{align*}
  q(q(w(t),\boldsymbol{y}), \boldsymbol{z}) & =_{\text{(ii)}} q(w(q(t_0,\boldsymbol{y}),\dots, q(t_k,\boldsymbol{y}),\dots),\boldsymbol{z}) \\
  & =_{\text{(ii)}} w(q (q (t_0,\boldsymbol{y}),\boldsymbol{z}),\dots,q (q(t_k,\boldsymbol{y}),\boldsymbol{z}),\dots)\\
  & =_{\text{ind.}} w(q(t_0, q(y_0, \boldsymbol{z}), \ldots),\ldots,q(t_k, q(y_k, \boldsymbol{z}), \ldots), \ldots ) \\
  & =_{\text{(ii)}} q (w(t),q (y_0,\boldsymbol{z}),\dots,q (y_k,\boldsymbol{z}),\dots) \text{.}
 \end{align*}
  Let $\mathcal C$ be an infinitary clone $\tau$-algebra and $\alpha:\dot{\mathscr{X}} \to C$ be a function. 
  We prove that there exists a unique homomorphism $\hat{\alpha} :\mathcal{N}_{\bar \tau}(\mathscr{X}) \to \mathcal C$ extending $\alpha$. 
  We define $\hat{\alpha}$ by induction:
  \begin{enumerate}
    \item[(a)] $\hat{\alpha}(\e_i^\mathcal{N})=\e_i^\mathcal C$;
    \item[(b)] $\hat{\alpha}(\mathscr{x}(t_0,\dots,t_k,\dots)) =  q^\mathcal C(\alpha(\dot{\mathscr{x}}),\hat{\alpha}(t_0),\dots,\hat{\alpha}(t_k),\dots)$ if $\mathscr{x} \in  \mathscr{X}$; 
    \item[(c)] $\hat{\alpha}(f(t_0,\dots,t_k,\dots))= q^\mathcal C(f^{\mathcal{C}},\hat{\alpha}(t_0),\dots,\hat{\alpha}(t_k),\dots)$ if $f \in \tau$.  
  \end{enumerate}
  The map $\hat{\alpha}$ is well-defined.
  We prove by induction that $\hat{\alpha}$ is a homomorphism:
  \begin{align*}
    \hat{\alpha}(q^\mathcal{N}(\mathscr{x}(t),\boldsymbol{u})) 
    & =_{\text{(ii)}}  \hat{\alpha}(\mathscr{x}(q^\mathcal{N}(t_0,\boldsymbol{u}),\dots, q^\mathcal{N}(t_k,\boldsymbol{u}),\dots)) \\
    & =_{\text{(b)}} q^\mathcal C(\alpha(\dot{\mathscr{x}}),\hat{\alpha}(q^\mathcal{N}(t_0,\boldsymbol{u})),\dots,\hat{\alpha}(q^\mathcal{N}(t_k,\boldsymbol{u})),\dots) \\
    & =_{\text{(ind.)}} q^\mathcal C(\alpha(\dot{\mathscr{x}}),q^\mathcal C(\hat{\alpha}(t_0),\hat{\alpha}^{\omega}(\boldsymbol{u})),\dots,q^\mathcal C(\hat{\alpha}(t_k),\hat{\alpha}^{\omega}(\boldsymbol{u})),\dots) \\
    & =_{\text{(N3)}} q^\mathcal C(q^\mathcal C(\alpha(\dot{\mathscr{x}}),\hat{\alpha}(t_0),\dots,\hat{\alpha}(t_k),\dots),\hat{\alpha}^{\omega}(\boldsymbol{u})) \\
    & =_{\text{(b)}} q^\mathcal C(\hat{\alpha}(\mathscr{x}(t)),\hat{\alpha}(u_0),\dots,\hat{\alpha}(u_k),\dots) \text{,}
  \end{align*}
  and similarly for $f \in \tau$. Let $\beta:\mathcal{N}_{\bar \tau}(\mathscr{X}) \to \mathcal C$ be another homomorphism extending $\alpha$. 
  Then we have, by induction
  \begin{align*}
    \beta(\mathscr{x}(t_0,\dots,t_k,\dots)) & =_{\text{(ii), (N2)}} \beta(q^\mathcal{N}(\dot{\mathscr{x}},t_0,\dots,t_k,\dots)) \\
    & =_{(\beta \text{ hom.})} q^\mathcal C(\beta(\dot{\mathscr{x}}),\beta(t_0),\dots,\beta(t_k),\dots ) \\
    & =_{\text{(ind.)}}q^\mathcal C(\alpha(\dot{\mathscr{x}}),\hat{\alpha}(t_0),\dots,\hat{\alpha}(t_k),\dots ) \\
    & =_{\text{(b)}} \hat{\alpha}(\mathscr{x}(t_0,\dots,t_k,\dots)) 
  \end{align*}
  for $\mathscr{x} \in  \mathscr{X}$. 
  Similarly for $f \in \tau$. 
\end{proof}
\end{proposition}

  \begin{remark}
    One could be tempted to define directly mutual translations between the two presentations of the free infinitary clone $\tau$-algebras: as equivalence classes of $\bar{\tau}$-terms or as $\tau$-metaterms.
One would then naturally define the translation of $q(v,\boldsymbol{t})$ by cases on $v$, and say in the case where $v$ is $q(\boldsymbol{s})$ that the translation
of $q(q(\boldsymbol{s}),\boldsymbol{t})$ is, by definition, the translation of $q(s_0,q(s_1,\boldsymbol{t}),\ldots,q(s_n,\boldsymbol{t}),\ldots)$, thus considering (N3) as a rewriting rule. 
It would be interesting  to explore the termination of this rewriting system. 
As a result of this discussion, we like to see our $\tau$-metaterms as $\bar{\tau}$-terms in ``normal form".
  \end{remark}

By Definition \ref{def:en}, there is a bijection between $N_{\bar \tau}$ and the set $T_{\tau}(X)$ of $\tau$-terms over a countable set of variables $X=\{x_0, x_1, \ldots\}$, 
that maps $\e_i$ to $x_i$ 
(see also Example \ref{exa:free}). 
As a consequence, it is easy to obtain the following.  

\begin{lemma}
  \label{lem_abs}
  The $\tau$-algebra $\mathcal{N}_{\bar \tau}^{\downarrow}$ is the absolutely free $\tau$-algebra over the set of generators $\{\e_0, \ldots, \e_n, \ldots \}$.
\end{lemma}

To each algebra we can canonically associate an infinitary clone algebra, the analogue of the clone of term operations. 

\begin{definition}
  \label{def:termclo}
Let $\mathbf{A}$ be a $\tau$-algebra. 
The \emph{term infinitary clone $\tau$-algebra $\mathbf{A}^{\uparrow}$}  
is the image of the unique  homomorphism $(-)^\mathbf{A}:\mathcal{N}_{\bar \tau}\to \mathcal{O}^{(\omega)}_\mathbf A$, 
that is the minimal subalgerba of $\mathcal{O}^{(\omega)}_\mathbf A$. 
\end{definition}

Explicitly, the term operation $t^\mathbf{A}$ is defined by induction as follows, for every $s\in A^\omega$:
\begin{equation*}
  t^\mathbf A(s)=\begin{cases}
    s^{\mathbf{A}}_i &\text{ if } t = \e_i\\
    f^\mathbf{A}(t_0^\mathbf{A}(s),\dots,t_n^\mathbf{A}(s),t_{n+1}^\mathbf{A}(s),\dots) & \text{ if } t = f(t_0,\dots,t_n,t_{n+1},\dots)
  \end{cases}
\end{equation*}

If $\mathsf{K}$ is a class of $\tau$-algebras, we denote by $\mathsf{K}^{\uparrow}$ the class $\{\mathbf{A}^{\uparrow} : \mathbf{A} \in \mathsf{K}\}$. 
We remark that the $\tau$-algebra $\mathbf A^{\updownarrows}$ is not in general isomorphic to $\mathbf A$.
Indeed, $\mathbf{A}^{\updownarrows}$ is always infinite but $\mathbf{A}$ may be finite. 

The kernel of the unique  homomorphism $(-)^\mathbf A:\mathcal{N}_{\bar \tau}\to \mathbf A^\uparrow$
coincides with the equational theory $\Th(\mathbf{A})$ of $\mathbf{A}$.
This holds by virtue of the identification between closed $\tau$-metaterms and $\tau$-terms expressed in Lemma \ref{lem_abs}. 

\begin{proposition}
  \label{prop:thm} 
  If $\mathcal C$ is a minimal infinitary clone $\tau$-algebra, then $\mathcal C^{\downuparrows} \simeq \mathcal C$. 
  \begin{proof}
    By definition,
    $\mathcal{C}^{\downuparrows}$ is the minimal subalgebra of $\mathcal{O}_{\mathcal{C}^{\downarrow}}^{(\omega)}$.  
    Since, by Proposition \ref{lem_repr}, $\mathcal{C}$ is embeddable into $\mathcal{O}_{\mathcal{C}^{\downarrow}}^{(\omega)}$, 
    then $\mathcal{C} \simeq \mathcal{C}^{\downuparrows}$. 
\end{proof}
\end{proposition}

As a consequence, for every congruence $\theta$ of $ \mathcal{N}_{\bar \tau}$, $(\mathcal{N}_{\bar \tau}/\theta)^{\downuparrows} \simeq \mathcal{N}_{\bar \tau}/\theta$ and for every $\tau$-algebra $\mathbf A$, $\mathbf A^{\updownarrows \uparrow} \simeq \mathbf A^\uparrow$. 

\section{Up-and-down between clone algebras and algebras}
\label{sec:up-down}
Let $\mathsf{K}$ be a class of $\tau$-algebras and $\mathsf{H}$ be a class of infinitary clone $\tau$-algebras. 
We define 
\begin{itemize}
  \item $\mathsf{K}^\uptriangle:= \mathbb I \mathbb S \{\mathcal{O}^{(\omega)}_{\mathbf{A}} : \mathbf{A} \in \mathsf{K}\}$; 
  \item $\mathsf{H}^{\downtriangle}:=\mathbb I \{\mathbf A: \mathsf{FCA}(\mathbf A) \cap \mathsf{H} \neq\emptyset  \}$.
\end{itemize}
We prove that 
\begin{itemize}
  \item if $\mathsf{K}$ is a variety of $\tau$-algebras, then $\mathsf{K}^{\uptriangle}$ is a variety of infinitary clone $\tau$-algebras;
  \item if $\mathsf{H}$ is a variety of infinitary clone $\tau$-algebras, then $\mathsf{H}^{\downtriangle}$ is a variety of $\tau$-algebras. 
\end{itemize}

The following lemma shows that the class of functional clone algebras is closed under products. 
\begin{lemma}
  \label{lem_prod}
  Let $\mathcal{F}_i \in \mathsf{FCA}(\mathbf{A}_i)$ and $\mathbf{A}:=\prod_{i \in I}\mathbf{A}_i$, for $i \in I$. 
  Then $\mathcal{F}:=\prod_{i \in I}\mathcal{F}_i$ can be embedded into $\mathcal{O}_\mathbf{A}^{(\omega)}$. 
  \begin{proof}
    Let $F=\prod_{i \in I}F_i$ and $A=\prod_{i \in I}A_i$. 
      Consider the function $\alpha$ that maps $a=(a(i): i \in I) \in F$ to $\alpha(a): A^\omega \to A$, where 
      \begin{equation*}
          \alpha(a)(b_0, b_1, \ldots, b_n, \ldots)(i)=a(i)(b_0(i), b_1(i), \ldots, b_n(i), \ldots) \text{.}
      \end{equation*}
      It is easy to see that $\alpha$ is an embedding. 
  \end{proof}
\end{lemma}

\begin{theorem} 
  \label{thm:ggg}
  If $\mathsf{K}$ is a variety of $\tau$-algebras, then $\mathsf{K}^\uptriangle$ is a variety of infinitary clone $\tau$-algebras. 
  \begin{proof}
      The fact that the class $\mathsf{K}^\uptriangle$ is closed under subalgebras is immediate. 
      Moreover, by Lemma \ref{lem_prod} and the hypothesis that $\mathsf{K}$ is a variety, the class $\mathsf{K}^\uptriangle$ is closed under products.
      We now prove that the class $\mathsf{K}^\uptriangle$ is closed under homomorphic images: let $\alpha: \mathcal{E} \to \mathcal{G}$ be a surjective homomorphism of infinitary clone algebras with $\mathcal{E} \in \mathsf{K}^{\uptriangle}$.
      Since $\mathcal{E} \in \mathsf{K}^{\uptriangle}$, there is a subalgebra $\mathcal{F}$ of $\mathcal{O}^{(\omega)}_{\mathbf{A}}$, for some $\mathbf{A} \in \mathsf{K}$, such that $\mathcal{E} \simeq \mathcal{F}$. 
      By Lemma \ref{rem_homisalso}, $\alpha$ is also a surjective homomorphism $\mathcal{E}^{\downarrow} \to \mathcal{G}^{\downarrow}$ and the isomorphism $\mathcal{E} \simeq \mathcal{F}$ is also an isomorphism $\mathcal{E}^{\downarrow} \simeq \mathcal{F}^{\downarrow}$. 
      By Proposition \ref{prop_rel}, $\mathcal{F}^{\downarrow}$ is embeddable into $\mathbf{A}^{A^\omega} \in \mathsf{K}$, so that $\mathcal{E}^{\downarrow}, \mathcal{F}^{\downarrow} \in \mathsf{K}$.
      Since $\mathsf{K}$ is a variety, $\mathcal{G}^{\downarrow} \in \mathsf{K}$ as well,
      so that $\mathcal{O}^{(\omega)}_{\mathcal{G}^{\downarrow}} \in \mathsf{K}^{\uptriangle}$. 
      By Proposition \ref{lem_repr} $\mathcal{G}$ can be embedded into $\mathcal{O}^{(\omega)}_{\mathcal{G}^{\downarrow}}$. 
      Hence $\mathcal{G} \in \mathsf{K}^{\uptriangle}$.
  \end{proof}
\end{theorem}

The following technical lemma is used in the proof of Theorem \ref{thm:giu}. 

\begin{lemma} 
  \label{lem:11}    
  Let $\alpha: \mathbf A\to \mathbf B$ be an onto homomorphism of $\tau$-algebras. The following conditions hold:
  \begin{enumerate}
    \item For every $\phi \in O^{(\omega)}_A$, if there is $\psi \in O^{(\omega)}_B$ such that $\psi\circ \alpha ^\omega= \alpha\circ \phi$, 
    then $\psi$ is unique, and is denoted by $\alpha^\star(\phi)$. 
    \item The set $C$ where $\alpha^{\star}$ is defined
    is a subalgebra of $\mathcal O^{(\omega)}_\mathbf A$.
    \item The map $\alpha^\star: C\to O^{(\omega)}_B$
      is a homomorphism from the functional clone $\tau$-algebra $\mathcal C$ of universe $C$ onto $\mathcal O^{(\omega)}_\mathbf B$. 
  \end{enumerate}

  \begin{proof}  
      (1) By surjectivity of $\alpha$. 

      \noindent
      (2) By hypothesis $\alpha^\omega : A^\omega\to B^\omega$ is onto. 
      Let $\phi_0,\phi_1,\dots,\phi_n,\ldots\in C$ and let $\psi_0,\psi_1,\dots,\psi_n,\ldots: B^\omega\to B$ such that $\psi_i\circ \alpha^\omega=\alpha\circ \phi_i$.
      For every $s \in A^{\omega}$: 
      \begin{align*}
        \alpha(q^A(\phi_0,\phi_1,\dots,\phi_n,\dots)(s)) & = \alpha(\phi_0(\phi_1(s),\dots,\phi_n(s),\dots))\\
        & = \psi_0( \alpha^\omega(\phi_1(s),\dots,\phi_n(s),\dots))\\
        & = \psi_0(\alpha(\phi_1(s)),\dots,\alpha(\phi_n(s)),\dots) \\
        & = \psi_0(\psi_1(\alpha^\omega(s)),\dots,\psi_n(\alpha^\omega(s)), \ldots)\\
        & = q^B(\psi_0,\psi_1,\dots,\psi_n,\dots)(\alpha^\omega(s)) \text{.} 
      \end{align*} 
      \noindent
      (3) The map $\alpha^\star$ can be proved to be a homomorphism by substituting $\alpha^\star(\phi_i)$ for $\psi_i$ in the chain of equalities of (2).
      We now prove that $\alpha^\star$ is onto. 
      Let $\psi \in O^{(\omega)}_B$. 
      For each $s \in A^\omega$, choose $\phi(s) \in \alpha^{-1}[\psi(\alpha^{\omega}(s))]$. 
      Then $\phi$ is the desired operation in $C$. \qedhere
\end{proof}
\end{lemma}

\begin{theorem}
  \label{thm:giu}
  If $\mathsf{H}$ is a variety of infinitary clone $\tau$-algebras, then $\mathsf{H}^{\downtriangle}$ is a variety of $\tau$-algebras.
 
  \begin{proof} (1) \emph{The class $\mathsf{H}^{\downtriangle}$ is closed under homomorphic images.}
  Let $\mathbf A\in  \mathsf{H}^{\downtriangle}$ and $\alpha: \mathbf A \to \mathbf B$ be an onto homomorphism. 
  Since $\mathbf A\in  \mathsf{H}^{\downtriangle}$, there exists $\mathcal{D} \in \mathsf{FCA}(\mathbf{A}) \cap \mathsf{H}$.
  Now, $\alpha$ is surjective, and by Lemma \ref{lem:11} there exist $\mathcal C \in \mathsf{FCA}(\mathbf{A})$ and a surjective homomorphism $\alpha^\star$ from $\mathcal C$ onto $\mathcal{O}^{(\omega)}_\mathbf B$.  
  Since $\mathcal D\in \mathsf{H}$, we have $\mathcal C\cap \mathcal D\in \mathsf{H}$. 
  Therefore $\alpha^\star[\mathcal C\cap \mathcal D] \in \mathsf{FCA}(\mathbf B) \cap \mathsf{H}$. 
  By definition of $\mathsf{H}^{\downtriangle}$ this implies that $\mathbf B\in  \mathsf{H}^{\downtriangle}$.

  (2) \emph{The class $ \mathsf{H}^{\downtriangle}$ is closed under subalgebras.}
  If $\mathbf A\in  \mathsf{H}^{\downtriangle}$, then there exists $\mathcal D \in \mathsf{FCA}(\mathbf{A}) \cap \mathsf{H}$. 
  Let $\mathbf B$ be a subalgebra of $\mathbf A$. 
  We consider the subalgebra $\mathcal C$ of  $\mathcal D$ of all $\phi \in D$ such that $\phi(s)\in B$ for every $s\in B^\omega$. 
  By hypothesis $\mathcal C\in \mathsf{H}$. 
  The image of the  homomorphism from $\mathcal C$ into $\mathcal O_\mathbf B^{(\omega)}$, defined by $\phi \mapsto \phi |_{B^\omega}$, belongs to $\mathsf{H}$ and is a subalgebra of $\mathcal O_\mathbf B^{(\omega)}$. 
  Consequently, $\mathbf B\in  \mathsf{H}^{\downtriangle}$.
 
  (3) \emph{The class $\mathsf{H}^{\downtriangle}$ is closed under products.}
  Let $\mathbf A=\prod_{i\in I}\mathbf A_i$ be the product of the algebras $\mathbf A_i\in  \mathsf{H}^{\downtriangle}$ and let $\mathcal C_i\in  \mathsf{H} \cap \mathsf{FCA}(\mathbf A_i)$ ($i\in I$). 
  By Lemma \ref{lem_prod}  $\prod_{i\in I} \mathcal C_i$ can be embedded into $\mathcal O_\mathbf{A}^{(\omega)}$. 
  Therefore $\mathbf A \in \mathsf{H}^{\downtriangle}$.  
\end{proof}
\end{theorem}

\section{A Birkhoff Theorem through infinitary clone algebras}
\label{sec:gb}
The classical Birkhoff Theorem for finitary algebras characterises the class $\mathsf{K}$ of models of an equational theory
(see e.g. \cite[Theorem 11.9]{BS81}); a proof for infinitary algebras is given in \cite[Theorem 9.6]{slo}.
This section contains an enhanced Birkhoff Theorem for algebras over a homogeneous type. 

\subsection{Closure under expansion}
To begin with, we prove 
that any $\tau$-algebra $\mathbf{A}$ can be embedded into a reduced product of its countably generated subalgebras. 
As a consequence we obtain that a variety of $\tau$-algebras is closed under expansion. 

\begin{definition}\label{def:As} Let $\mathbf A$ be an $\tau$-algebra.
The \emph{subalgebra $\mathbf A_{\bar s}$} of $\mathbf A$ \emph{generated} by $s\in A^\omega$ is  the image of the homomorphism $ \bar{s}: t^\mathbf A \mapsto t^\mathbf{A}(s)$ (Proposition \ref{prop_rel}(1)), that 
is the subalgebra of $\mathbf{A}$ generated by $\set(s)$. 
\end{definition}

The subalgebra $\mathbf A_{\bar s}$ is finitely generated if $|\set(s)|<\omega$; 
otherwise, it has a countable infinite set of generators. 

\begin{definition} Let $\mathsf{K}$ be a class of $\tau$-algebras.
We say that 
\begin{enumerate}
  \item $\mathsf{K}$ is \emph{closed under expansion} (resp. \emph{under finite expansion}) 
if, for every $\tau$-algebra $\mathbf A$,
\begin{equation*}
  (\forall s\in A^\omega \, \mathbf A_{\bar s} \in \mathsf{K}) \IMP \mathbf A\in \mathsf{K} \quad \text{ (resp. } (\forall s\in A_{\fin}^\omega \, \mathbf A_{\bar s} \in \mathsf{K}) \IMP \mathbf A\in \mathsf{K})  \text{.}
\end{equation*}
\end{enumerate}
\end{definition}

\begin{theorem}
  \label{thm_exp}
  If $\mathsf{K}$ is a variety of $\tau$-algebras, then $\mathsf{K}$ is closed under expansion. 
  \begin{proof}
      Let $\mathbf{A}$ be a  $\tau$-algebra.  
      Assume that for every $s \in A^\omega$, $\mathbf A_{\bar s}\in \mathsf{K}$. 
      Let $I:=A^{\omega}$. 
      For $s \in I$, let $J_s:=\{i \in I : \set(s) \subseteq \set(i)\}$. 
      The family $\{J_s : s \in A^{\omega}\}$ is closed under countable intersections: if $s_0, s_1, s_2, \ldots \in A^{\omega}$,
      then $\bigcap\{ J_{s_k} : k \in \omega \} = J_r$
      where $r \in A^{\omega}$ is any sequence such that $\set(r)=\bigcup\{\set(s_k) : k \in \omega\}$.
      Let $\EuScript{F}$ be the closure of the family $\{J_s : s \in A^{\omega}\}$ by supersets. 
      Then $\EuScript{F}$ is a proper countably complete filter on $I$. 
      Observe that, for $x,y \in B:=\prod_{s \in I} A_{\bar{s}}$ the relation 
       $ (x,y) \in \theta_{\EuScript{F}} \IFF \{i \in I : x(i) = y(i)\} \in \EuScript{F}$
      is a congruence. 
      For $a \in A$, let $\lambda a$ be any element of $B$ such that $(\lambda a)(s)=a$ if $a \in \set(s)$.
      We show that the function $\alpha: A \to B / \theta_{\EuScript{F}}$, 
      defined by $\alpha(a)=(\lambda a)/\theta_{\EuScript{F}}$, is an embedding. 
      
      (1) \emph{The map $\alpha$ is well-defined}. 
      If $\lambda' a$, with $(\lambda'a)(i) = a $ if $a \in \set(i)$, is used instead of $\lambda a$ then 
      \begin{equation*}
          J_{(a,a,\ldots)} = \{i \in I : a \in \set(i)\} \subseteq \{i \in I : (\lambda a)(i) = (\lambda'a)(i)\} \in \EuScript{F}
      \end{equation*}
      so that $(\lambda'a, \lambda a) \in \theta_{\EuScript{F}}$.

      (2) \emph{The map $\alpha$ is a homomorphism}. 
      For all $s=(a_0, a_1, \ldots, a_k, \ldots) \in A^\omega$ and $f \in \tau$, 
      \begin{equation*}
          f^{\mathbf{B} / \theta_{\EuScript{F}}}(\lambda a_0/\theta_{\EuScript{F}},\ldots, \lambda a_k /\theta_{\EuScript{F}}, \ldots) 
          = f^\mathbf{B}(\lambda a_0, \ldots, \lambda a_k, \ldots) / \theta_{\EuScript{F}} 
      \end{equation*}
      and 
          $\alpha(f^\mathbf{A}(s)) = \lambda f^\mathbf{A}(s) / \theta_{\EuScript{F}} \text{,}$
      so that it is enough to show that $\lambda f^\mathbf{A}(s)$ and $f^\mathbf{B}(\lambda a_0, \ldots, \lambda a_k, \ldots)$ are congruent modulo $\theta_{\EuScript{F}}$: 
      \begin{align*}
          J_s & = \{i \in I: \set(s) \subseteq \set(i)\} \\
          & \subseteq \{i \in I : f^{\mathbf{A}_{\bar{i}}}(s) =f^\mathbf{A}(s) \} \\
          & = \{i \in I : f^{\mathbf{A}_{\bar{i}}}(\lambda a_0(i), \ldots, \lambda a_k(i), \ldots) =\lambda f^\mathbf{A}(s)(i)\}\\
          & = \{i \in I: f^\mathbf{B}(\lambda a_0, \ldots, \lambda a_k, \ldots)(i) = \lambda f^\mathbf{A}(s)(i)\} \text{.}
      \end{align*}
      (3) \emph{The map $\alpha$ is injective}. Let $a,b$ be distinct elements of $A$; we need to show that $\{i \in I: (\lambda a)(i) = (\lambda b)(i)\} \notin \EuScript{F}$. 
      For every $i \in I$ such that $\{a,b\} \subseteq \set(i)$, $(\lambda a)(i)=a \neq b=(\lambda b)(i)$.
      Since $J_{(a,b,a,b, \ldots)} \in \EuScript{F}$ and 
      \begin{equation*}
        J_{(a,b,a,b, \ldots)} = \{i \in I: a,b  \in \set(i)\} \subseteq \{i \in I : (\lambda a)(i) \neq (\lambda b)(i)\} 
      \end{equation*}
      we conclude that $\{i \in I : (\lambda a)(i) \neq (\lambda b)(i)\} \in \EuScript{F}$.
      If, towards a contradiction, $\{i \in I: (\lambda a)(i) = (\lambda b)(i)\} \in \EuScript{F}$, then $\emptyset \in \EuScript{F}$. 
      Absurd. \qedhere
  \end{proof}
\end{theorem}

\begin{definition}
  Let $\mathbf{A}$ be a $\tau$-algebra. 
  We say that $\mathbf{A}$ is \emph{finite-dimensional} if $\mathrm{Fd}(\mathbf{A}^{\uparrow})=\mathbf{A}^{\uparrow}$ (see Definition \ref{def:dim}). 
\end{definition}

The proof of Theorem \ref{thm_exp} can be adapted to obtain the following. 

\begin{proposition}
  \label{prop_exp}
  Let $\mathsf{K}$ be a variety of finite-dimensional $\tau$-algebras. 
  Then $\mathsf{K}$ is closed under finite expansion. 
  \begin{proof}
    Let $\mathbf{A}$  be a finite-dimensional $\tau$-algebra.  
    Assume that for every $s \in A_{\fin}^\omega$, $\mathbf A_{\bar s}\in \mathsf{K}$. 
    Let $I:=A_{\fin}^{\omega}$. 
    For $s \in I$, let $J_s:=\{i \in I : \set(s) \subseteq \set(i)\}$. 
    The family $\{J_s : s \in A^{\omega}\}$ is closed under finite intersections,
    hence if 
    $\EuScript{F}$ is the closure of the family $\{J_s : s \in A^{\omega}\}$ by supersets, 
    $\EuScript{F}$ is a proper filter on $I$. 
    As $\mathbf{A}$ is finite-dimensional, it is not necessary that $\EuScript{F}$ is countably-complete to repeat the argument given above. 
  \end{proof}
\end{proposition}

Secondly, we recall the well-known result that equational classes are varieties,
with a proof homogeneous with our presentation. 

\begin{lemma}\label{lem:ti} Let $\Sigma$ be a set of $\tau$-identities between $\tau$-terms. 
  Then $\Mod(\Sigma)=\{\mathbf{A} \in \mathsf{Alg}_{\tau}: \Sigma\subseteq \Th(\mathbf A)\}$ is a variety.
  
    \begin{proof}
      (1) \emph{The class $\Mod(\Sigma)$ is closed under subalgebras.}
      Let $\mathbf{B}$ be a subalgebra of $\mathbf{A} \in \Mod(\Sigma)$. 
      The restriction function $(-|B): \mathbf{A}^{\uparrow} \to \mathbf{B}^{\uparrow}$
      is a homomorphism of infinitary clone $\tau$-algebras. 
      Now, since $\mathcal{N}_{\bar \tau}$ is initial we have the following factorisation 
        \begin{equation*}
          \begin{tikzcd}
              \mathcal{N}_{\bar \tau} \arrow[rr, "(-)^{\mathbf{B}}"] \arrow[dr, "(-)^{\mathbf{A}}", swap] & &  \mathbf{B}^{\uparrow}\\
              & \mathbf{A}^{\uparrow} \arrow[ur, "(-|B)", swap] & \\
          \end{tikzcd}
        \end{equation*}
        \vspace{-0.75cm}

      \noindent
      Observing that $\ker ((-)^\mathbf{A}) \subseteq \ker ((-|B) \circ (-)^\mathbf{A}) $, the above factorisation implies that $\Th(\mathbf{A}) \subseteq \Th(\mathbf{B})$. 
      Since $\Sigma \subseteq \Th(\mathbf{A})$, $\Sigma \subseteq \Th(\mathbf{B})$. 

      (2) \emph{The class $\Mod(\Sigma)$ is closed under products.}
      Let $\mathbf{A}_i \in \Mod(\Sigma)$ for every $i \in I$ and $\mathbf{A}:=\prod_{i \in I}\mathbf{A}_i$. 
      By Lemma \ref{lem_prod}, $\prod_{i \in I} \mathbf{A}^{\uparrow}_i$ is (isomorphic to) a subalgebra of $\mathcal{O}_{\mathbf{A}}^{(\omega)}$. 
      Since $\mathcal{N}_{\bar \tau}$ is initial the homomorphism
      \begin{equation*}
          ((-)^{\mathbf{A}_i}:i \in I) : \mathcal{N}_{\bar \tau} \to \prod_{i \in I} \mathbf{A}^{\uparrow}_i
      \end{equation*}
      factors as 
      \begin{equation*}
        \begin{tikzcd}
            \mathcal{N}_{\bar \tau} \arrow[rr] \arrow[dr, "(-)^{\mathbf{A}}", swap] & &  \prod_{i \in I} \mathbf{A}^{\uparrow}_i\\
            & \mathbf{A}^{\uparrow} \arrow[ur, hook] & \\
        \end{tikzcd}
      \end{equation*}
      \vspace{-0.75cm}

      \noindent
      Observing that $\ker((-)^{\mathbf{A}_i}:i \in I) = \bigcap_{i \in I} \Th(\mathbf{A}_i)$,
      the above diagram implies that $\Th(\mathbf{A}) = \bigcap_{i \in I} \Th(\mathbf{A}_i) \text{.}$
      Since $\Sigma \subseteq \Th(\mathbf{A}_i)$ for all $i \in I$, $\Sigma \subseteq \Th(\mathbf{A})$. 

      (3) \emph{The class $\Mod(\Sigma)$ is closed under homomorphic images.}
      Let $\alpha: \mathbf{A} \to \mathbf{B}$ be an onto homomorphism with $\mathbf{A} \in \Mod(\Sigma)$.
      Let $\mathcal{C}$ and $\alpha$ as in Lemma \ref{lem:11}; 
      $\mathbf{A}^{\uparrow}$ is a subalgebra of $\mathcal{C}$, and  
      by induction on $t \in \mathcal{N}_{\bar \tau}$, we see that $t^{\mathbf{A}} \circ \alpha^{\omega} = \alpha \circ t^{\mathbf{B}}$, i.e. $\alpha^{\star}(t^{\mathbf{A}}) = t^{\mathbf{B}}$. 
      Therefore, the homomorphism $\alpha^{\star} : \mathcal{C} \to \mathcal{O}_{\mathbf{B}}^{(\omega)}$ of Lemma \ref{lem:11} restricts to a homomorphism 
      $\mathbf{A}^{\uparrow} \to \mathbf{B}^{\uparrow}$.
      With the same argument given in the proof of (1), we conclude that $\Sigma \subseteq \Th(\mathbf{B})$. \qedhere
  \end{proof}
  
\end{lemma}

\subsection{An enhanced Birkhoff Theorem}
The equivalence of (1) and (2) in Theorem \ref{thm:bir1} below is Birkhoff Theorem (for algebras over a homogeneous type).
The third item (3) expands this equivalence with a new condition, and allows for a new proof.  

\begin{theorem}\label{thm:bir1} 
  Let $\mathsf{K}$ be a class of $\tau$-algebras. 
  Then the following conditions are equivalent:
    \begin{enumerate}
      \item $\mathsf{K}$ is a variety.
      \item $\mathsf{K}=\Mod(\Th(\mathsf{K}))$.
      \item $\mathsf{K}=\mathsf{K}^{\uptriangle \downtriangle}$ and $\mathsf{K}^\uptriangle$ is a variety of infinitary clone $\tau$-algebras.
 \end{enumerate}


  \begin{proof}
  $(1) \IMP (2)$ If $\mathbf{A} \in \mathsf{K}$, then by definition $\Th(\mathsf{K}) \subseteq \Th(\mathbf{A})$, so that $\mathbf{A} \in \Mod(\Th(\mathsf{K}))$. 
  For the opposite inclusion, let $\mathcal{N}_{\mathsf{K}} := \mathcal{N}_{\bar \tau} / \Th(\mathsf{K})$. 
  Firstly, we prove that $\mathcal{N}_{\mathsf{K}}^\downarrow \in  \mathsf{K}$. 
  Let $I$ be the set of all $\tau$-identities $i$ such that $i$ fails in some $\mathbf A_i\in \mathsf{K}$. 
  Now, we prove that $\Th(\mathsf{K})$ is the kernel $\theta_I$ of the unique homomorphism $\mathcal{N}_{\bar \tau} \to \prod_{i \in I} \mathbf{A}^{\uparrow}_i$. 
  Obviously, $\Th(\mathsf{K}) \subseteq \theta_I$ and the inclusion cannot be proper as by definition of $I$ this would lead to a contradiction. 
  Hence, $\mathcal{N}_{\mathsf{K}}$ is a subalgebra of $\prod_{i \in I} \mathbf{A}^{\uparrow}_i$. 
  Consider the product $\mathbf{B}:=\prod_{i \in I} \mathbf{A}_i \in \mathsf{K}$.  
  By Lemma \ref{lem_prod} the infinitary clone $\tau$-algebra $\prod_{i \in I} \mathbf{A}^{\uparrow}_i$ is isomorphic to some $\mathcal{F} \in \mathsf{FCA}(\mathbf{B})$.
  By Proposition \ref{prop_rel} $\mathcal{F}^{\downarrow}$ is embeddable in some power of $\mathbf{B}$, hence belongs to $\mathsf{K}$. 
  Since $\mathcal{N}_{\mathsf{K}}^{\downarrow}$ is a subalgebra of $\mathcal{F}^{\downarrow}$, $\mathcal{N}_{\mathsf{K}}^{\downarrow}$ belongs to $\mathsf{K}$ as well. 

  Now we prove that if $\mathbf{A} \in \Mod(\Th(\mathsf{K}))$, then $\mathbf{A} \in \mathsf{K}$. 
  Since $\Th(\mathsf{K}) \subseteq \Th(\mathbf{A})$, the unique surjective homomorphism $\mathcal{N}_{\bar \tau} \to \mathbf{A}^{\uparrow}$ can be factored through $\mathcal{N}_{\mathsf{K}}$. 
  The homomorphism $\mathcal{N}_{\mathsf{K}} \to \mathbf{A}^{\uparrow}$ is also surjective, as well as $\mathcal{N}_{\mathsf{K}}^{\downarrow} \to \mathbf{A}^{\updownarrows}$ by Lemma \ref{rem_homisalso}. 
  Since $\mathsf{K}$ is a variety, and since $\mathcal{N}_{\mathsf{K}}^{\downarrow} \in \mathsf{K}$, we get $\mathbf{A}^{\updownarrows} \in \mathsf{K}$. 
  Now, for every $s \in A^\omega$, by Proposition \ref{prop_rel}, $\bar{s} : \mathbf{A}^{\updownarrows} \to \mathbf{A}$ is a homomorphism, so that $\bar{s}[\mathbf{A}^{\updownarrows}]$, which is the subalgebra $\mathbf A_{\bar s}$ generated by $s$, belongs to $\mathsf{K}$. 
  Theorem \ref{thm_exp} implies that $\mathbf{A} \in \mathsf{K}$. 
  This concludes the proof that  $\mathsf{K}=\Mod(\Th(\mathsf{K}))$.

  \noindent
  (2) $\IMP$ (3) 
  If $\mathsf{K}=\Mod(\Th(\mathsf{K}))$,
  by Lemma \ref{lem:ti}, 
  $\mathsf{K}$ is a variety. 
  Then by Theorem \ref{thm:ggg} $\mathsf{K}^{\uptriangle}$ is a variety of infinitary clone $\tau$-algebras.
  By definition, for any class $\mathsf{K}$, we have $\mathsf{K}\subseteq \mathsf{K}^{\uptriangle \downtriangle}$. 
  We now show that $\mathsf{K}^{\uptriangle \downtriangle} \subseteq \mathsf{K}$. 
  Let $\mathbf{A} \in \mathsf{K}^{\uptriangle \downtriangle}$; 
  then there is $\mathcal{G} \in \mathsf{K}^\uptriangle \cap \mathsf{FCA}(\mathbf{A})$.
  Consequently, there is $\mathbf{B} \in \mathsf{K} = \Mod(\Th(\mathsf{K}))$, and 
  $\mathcal{F} \in \mathsf{FCA}(\mathbf{B})$ such that $\mathcal{G} \simeq \mathcal{F}$.
  This isomorphism restricts to 
  an isomorphism $\mathbf{A}^{\uparrow} \simeq \mathbf{B}^{\uparrow}$ of the minimal subalgebras, so that $\Th(\mathsf{K}) \subseteq \Th(\mathbf{B}) = \Th(\mathbf{A})$. 
  It follows that $\mathbf{A} \in \mathsf{K}$. 

  \noindent
  (3) $\IMP$ (1) By Theorem \ref{thm:giu}.
\end{proof}
\end{theorem}

As a byproduct, we obtain, for a variety $\mathsf{K}$, a different characterisation of the free $\mathsf{K}$-algebra over a countable set of generators.  
\begin{lemma}
  \label{lem_free}
  Let $\mathsf{K}$ be a variety of $\tau$-algebras with theory $\Th(\mathsf{K})$.
  Let $\mathcal{N}_{\mathsf{K}}:=\mathcal{N}_{\bar \tau}/\Th(\mathsf{K})$. 
  Then $\mathcal{N}_{\mathsf{K}}^{\downarrow}$ is the free $\mathsf{K}$-algebra over
  $\{\e_0, \ldots, \e_n, \ldots \}$.
  \begin{proof}
    The fact that $\mathcal{N}_{\mathsf{K}}^{\downarrow} \in \mathsf{K}$ has already been showed in the proof of (1) $\IMP$ (2), Theorem \ref{thm:bir1}. 
    Applying Lemma \ref{rem_homisalso} to the quotient map $\mathcal{N}_{\bar \tau} \to \mathcal{N}_{\mathsf{K}}$ gives a surjective 
    homomorphism  $\mathcal{N}_{\bar \tau}^{\downarrow} \to \mathcal{N}_{\mathsf{K}}^{\downarrow}$. 
    Consequently, $\mathcal{N}_{\bar \tau}^{\downarrow} /\Th(\mathsf{K}) \simeq \mathcal{N}_{\mathsf{K}}^{\downarrow}$ and we conclude thanks to Lemma \ref{lem_abs}. 
  \end{proof} 
\end{lemma}

Note that neither the first nor the second condition of the third item of Theorem \ref{thm:bir1} can be removed, as Examples \ref{ex_cond1} and \ref{ex_cond2} show.  

\begin{example}
  \label{ex_cond1}
  Let $\tau = \emptyset$ be the empty type. 
  The class $\mathsf{Set}$ of all sets is the variety of all algebras over $\tau$. 
  The class $\mathsf{Set}^{\uptriangle}$ is the class $\mathsf{CA}_{\emptyset}$ of all pure infinitary clone algebras. 
  Now, if  $A \in \mathsf{Set}$  and  $\mathsf{K}:=\mathsf{Set} \setminus \{A\}$,
  then $\mathsf{K}^{\uptriangle} = \mathsf{CA}_{\emptyset}$, but $\mathsf{K} \neq \mathsf{K}^{\uptriangle \downtriangle}$ is not a variety.   
  This shows that the hypothesis that $\mathsf{K}^{\uptriangle }$ is a variety alone does not imply that $\mathsf{K}$ is a variety.
\end{example}

If $\mathsf{K}$ is a class of $\tau$-algebras we shall denote by $\mathsf{K}^{\parallelogramm}$ the class $\mathsf{K}^{\uptriangle \downtriangle}$. 
If $\mu \in \mathsf{Ord}$ is an ordinal, we define $(-)^{\parallelogramm^\mu}$ as a class operator defined by (transfinite) induction: 
\begin{itemize}
  \item $\mathsf{K}^{\parallelogramm^0} =\mathsf{K}$; 
  \item $\mathsf{K}^{\parallelogramm^\mu} = (\mathsf{K}^{\parallelogramm^\nu})^{\parallelogramm}$ if $\mu = \nu +1$; 
  \item $ \mathsf{K}^{\parallelogramm^\mu} = \bigcup\{ \mathsf{K}^{\parallelogramm^\nu} : \nu < \mu\}$ if $\mu$ is limit. 
\end{itemize} 
By definition, we have that $\mathsf{K} \subseteq \mathsf{K}^{\parallelogramm}$, and if $\mathsf{K}$ is a variety, by Theorem \ref{thm:bir1}, equality holds. 
Iterating, we immediately get that $\mathsf{K} \subseteq \mathsf{K}^{\parallelogramm^\mu}$ for every ordinal $\mu$. 
One may wonder whether $\bigcup\{\mathsf{K}^{\parallelogramm^\mu} : \mu \in \mathsf{Ord} \}$ 
is the smallest variety containing $\mathsf{K}$.
This is not the case, as we show in Example \ref{ex_cond2}.   
\begin{lemma}
  Let $\mathsf{K}$ be a class of $\tau$-algebras. 
  Then 
  \begin{equation}
    \label{eq_paral}
    \mathsf{K} \subseteq \bigcup_{\mu \in \mathsf{Ord}} \mathsf{K}^{\parallelogramm^\mu}  \subseteq \bigcup_{\mathbf{A} \in \mathsf{K}}\mathbb{V}(\mathbf{A}) \subseteq \mathbb{V}(\mathsf{K}) \text{.}
  \end{equation}
  \begin{proof}
    The proof is
    by (transfinite) induction on $\mu$. 
    If $\mu=0$, the conclusion trivially follows. 
    Assume that \eqref{eq_paral} holds for every $\nu < \mu$. 
    Let $\mathbf{B} \in \mathsf{K}^{\parallelogramm^\mu}$. 
    Then there are $\nu < \mu$ and a subalgebra $\mathcal{F}$ of $\mathcal{O}_{\mathbf{B}}^{(\omega)}$ such that $\mathcal{F} \in (\mathsf{K}^{\parallelogramm^\nu})^{\uptriangle}$. 
    By definition there is $\mathbf{A}\in \mathsf{K}^{\parallelogramm^\nu}$ such that (up to isomorphism) $\mathcal{F}$ is a subalgebra of $\mathcal{O}_{\mathbf{A}}^{(\omega)}$. 
    This implies that $\mathbf{A}^{\uparrow} \simeq \mathbf{B}^{\uparrow}$ and, consequently, that $\Th(\mathbf{A}) = \Th(\mathbf{B})$. 
    Therefore, by inductive hypothesis, $\mathbf{B} \in \Mod(\Th(\mathbf{A}))=\mathbb{V}(\mathbf{A})$.   
  \end{proof}
\end{lemma}

\begin{example}
  \label{ex_cond2}
  Let $\mathbf{A}$ and $\mathbf{B}$ be two $\tau$-algebras for some homogeneous type $\tau$. 
  Let $\mathsf{K}:=\{\mathbf{A}, \mathbf{B}\}$ and $\mathsf{L}:=\bigcup\{\mathsf{K}^{\parallelogramm^\mu} : \mu \in \mathsf{Ord}\}$. 
  The class $\mathsf{L}$ is not in general a variety, as $\mathsf{L} \subseteq \mathbb{V}(\mathbf{A}) \cup \mathbb{V}(\mathbf{B}) \subset \mathbb{V}(\mathsf{K})$. 
  The second inclusion is strict since generally $\mathbf{A} \times \mathbf{B} \in \mathbb{V}(\mathsf{K}) \setminus (\mathbb{V}(\mathbf{A}) \cup \mathbb{V}(\mathbf{B}))$.  
  However, 
  $\mathsf{L}^{\parallelogramm}=\mathsf{L}$ but $\mathsf{L}$ is not a variety, 
  implying that the first condition of (3) of Theorem \ref{thm:bir1} alone is not enough to derive (1).  
  This also implies that $\bigcup\{\mathsf{K}^{\parallelogramm^\mu} : \mu \in \mathsf{Ord} \}$ 
  is not the smallest variety containing $\mathsf{K}$.
\end{example}

\subsection{Hyperidentities} 
\label{subs:central}
The purpose of this section is to explore the relationship between $\Th({\mathsf{K}^\uptriangle})$ and $\Th({\mathsf{K}^{\uparrow}})$ for a given class of algebras $\mathsf{K}$. 
Along the way, we will cast some light on the variety $\mathbb{V}(\mathcal{P})$ generated by the pure infinitary clone algebra $\mathcal P$ defined in Example \ref{exa:proj}.
If $\mathsf{K}$ is class of $\tau$-algebras, recall that $\mathsf{K}^{\uparrow}=\{\mathbf A^{\uparrow} : \mathbf A\in \mathsf{K}\}$. 
Since $\mathbf{A}^{\uparrow}$ is a subalgebra of $\mathcal{O}^{(\omega)}_{\mathbf{A}}$, $\mathsf{K}^{\uparrow}$ is a subclass of $\mathsf{K}^{\uptriangle}$. 
Then $\Th({\mathsf{K}^\uptriangle}) \subseteq \Th({\mathsf{K}^{\uparrow}})$.
We will show in Example \ref{ex:central} that in general $\Th({\mathsf{K}^\uptriangle})\neq \Th({\mathsf{K}^{\uparrow}})$. 

Let $\mathcal{C}$ be an infinitary clone $\tau$-algebra.
A sequence $(\theta_0, \theta_1, \ldots, \theta_n, \ldots )$ of congruences of $\mathcal{C}$ is a sequence of \emph{complementary factor congruences}
if $\bigcap\{\theta_i : i \in \omega\} = \Delta_{C}$ and for every $s \in C^{\omega}$ there is $b \in C$ such that $(b,s_i) \in \theta_i$. 
An element $a \in C$ is called \emph{central} if the principal congruences $\theta(a, \e_i)$, $i \in \omega$, generated by $(a, \e_i)$, form a sequence of complementary factor congruences. 
A central element $a$ is \emph{nontrivial} if $a \neq \e_i$ for all $i \in \omega$.

\begin{lemma}
  \label{lemma_central}
  Let $\mathcal{C}$ be an infinitary clone $\tau$-algebra. 
  An element $a \in C$ is central iff the following conditions hold:
  \begin{itemize}
    \item[(C1)] $q(a,x,\dots,x,\dots) = x$;
    \item[(C2)] $q(a,q(a,x_{0}^0,x_{0}^1,\ldots),q(a,x_{1}^0,x_{1}^1,\ldots),q(a,x_{2}^0,x_{2}^1,\ldots),\ldots)$ \\
      $ \text{} \hspace{1cm} =q(a,x_{0}^0,x_{1}^1,x_{2}^2,\ldots)$.
    \item[(C3)] $q(a, q(y_0, z_0^0, z_1^0, \ldots), q(y_1, z_0^1, z_1^1, \ldots), \ldots )$\\ 
      $ \text{} \hspace{1cm} =q(q(a, \boldsymbol{y}), q(a, z_0^0, z_0^1, \ldots), q(a, z_1^0, z_1^1, \ldots), \ldots )$.
  \end{itemize}
  \begin{proof}
    Easy generalisation of \cite[Proposition 3.6]{SLPK13}. 
  \end{proof}
\end{lemma}

If $\mathcal{C} \in \mathbb{V}(\mathcal{P})$, then every element of $\mathcal{C}$ is central. 
In particular $\mathbb{V}(\mathcal{P})$ satisfies the identities (C1)--(C3). 

\begin{example}
  \label{ex:central}
If $\mathsf{Set}$ is the variety of sets, then $\mathsf{Set}^{\uptriangle}$
is the class of pure infinitary clone algebras, 
while $\mathsf{Set}^{\uparrow}=\mathbb{I}(\mathcal{P}) \subseteq \mathbb{V}(\mathcal{P})$.
Since there are pure infinitary clone algebras in which not every element is central, 
the identities \emph{(C1)--(C3)}
belong to $\Th({\mathsf{Set}^{\uparrow}}) \setminus \Th({\mathsf{Set}^\uptriangle})$.
\end{example}

Given a class $\mathsf{K}$ of $\tau$-algebras, the identities satisfied by $\mathsf{K}^{\uptriangle}$ are very few: 
for an identity to be valid in $\mathsf{K}^{\uptriangle}$, it must hold true in $\mathcal{O}^{(\omega)}_{\mathbf{A}}$ for every $\mathbf{A} \in \mathsf{K}$.
In contrast, the identities satisfied by $\mathsf{K}^{\uparrow}$ are more interesting. 
These encompass precisely those identities satisfied by $\mathbf{A}^{\uparrow}$ for every $\mathbf{A} \in \mathsf{K}$, and correspond to the set of hyperidentities satisfied by $\mathsf{K}$.
Hyperidentities, as outlined in \cite{T81} and \cite[Section 10.10]{ALV3}, expand the notion of identities to include variables that range over term operations as well as elements of the algebras. 

\section{An enhanced Birkhoff Theorem for finitary algebras}
\label{sec:tafromalg}

\subsection{Homogenising and dehomogenising terms} 
In what follows $\rho=(\rho_n:  n \in \omega)$ is a finitary type.
We associate with $\rho$ the homogeneous type $\rho_\bullet:=\bigcup\{\rho_n : n \in \omega\}$.
Otherwise stated, if $f \in \rho_n$ is an operation symbol of arity $n$, $f \in \rho_{\bullet}$ is a symbol of arity $\omega$. 
Let $\mathbf S=(S,f^\mathbf S)_{f\in\rho}$ be a $\rho$-algebra. 
We denote by $\mathbf S^\top= (S, f^{(\mathbf{S}^\top)})_{f\in\rho_{\bullet}}$ the $\rho_\bullet$-algebra such that $f^{(\mathbf{S}^\top)}$ is the top extension of $f^\mathbf S$
(see Definition \ref{def:topext}).
If $\mathsf{H}$ is a class of $\rho$-algebras, we define $\mathsf{H}^{\top}:=\{\mathbf{S}^{\top} : \mathbf{S} \in \mathsf{H}\}$.

\begin{definition}
  \label{def:homo}
  Let $T_\rho(V)$ be  the set of $\rho$-terms over $V=\{\e_0,\e_1,\dots\}$. 
  We define two maps, $(-)^\bullet: T_\rho(V)\to N_{\bar {\rho_\bullet}}$ and $(-)^\circ: N_{\bar{\rho_\bullet}} \to T_\rho(V)$ as follows:
  \begin{enumerate}
    \item $(\e_i)^\bullet= \e_i$ and $(\e_i)^\circ= \e_i$ for every $i \in \omega$;
    \item $f(t_0,\dots,t_{n-1})^\bullet= f(t_0^\bullet,\dots,t_{n-1}^\bullet,\e_{n},\e_{n+1},\dots)$, for every $f \in\rho_n$;  
    \item $f(t_0,\dots,t_{n-1},t_{n},t_{n+1},\dots)^\circ=f(t_0^\circ,\dots,t_{n-1}^\circ)$, for every $f\in\rho_n$. 
  \end{enumerate}
\end{definition}

The map $(-)^\circ$ is onto, while $(-)^\bullet$ is not onto. If $\Sigma$ is a set of $\rho$-identities, we write $\Sigma^\bullet$ for the set $\{t^\bullet = u^\bullet : t=u \in \Sigma\}$. 

\begin{lemma}
  \label{lem:upterms} 
  Let $\mathbf S$ be a $\rho$-algebra. 
  Then, for every closed $\rho_\bullet$-metaterms $t,u$ and $\rho$-terms $v,z$:
  \begin{itemize}
    \item[(i)] $\mathbf S^\top \models (t^\circ)^\bullet  = t$;
    \item[(ii)] the term $(v^\bullet)^\circ$ is syntactically equal to $v$;
    \item[(iii)] $\mathbf S^\top\models t=u  \IFF \mathbf S^\top\models (t^\circ)^\bullet=(u^\circ)^\bullet \IFF \mathbf S \models t^\circ=u^\circ$;
    \item[(iv)] $\mathbf S^\top\models v^\bullet=z^\bullet \IFF  \mathbf S \models v=z$. 
\end{itemize}
\begin{proof} 
  The proof of (i) and (ii) is straightforward;
  (iii) follows from (i) and the definition of $\mathbf S^\top$;
  (iv) follows from (ii) and (iii) letting $v=t^\circ$ and $z = u^\circ$.
\end{proof}
\end{lemma}

\begin{definition} 
  Let $t\in N_{\bar{\rho_\bullet}}$. 
  The $\rho_\bullet$-identity $t= (t^\circ)^\bullet $  is called \emph{$\rho$-structural}. 
  The set of all $\rho$-structural identities associated with the finitary type $\rho$ will be denoted by $\Str(\rho)$.
\end{definition}

\begin{example}
  \label{exa:stru}
  Let $\rho$ be a finitary type. 
  If $f \in \rho_n$ for $n\geq 1$, the identity 
  \begin{equation}
    \label{eq:findim}
  f(\e_0,\ldots,\e_{n-1},\e_0,\e_0,\ldots)=f(\e_0,\e_1, \ldots, \e_{n-1}, \e_{n}, \e_{n+1},\ldots)
  \end{equation}
  specifying that $\dim(f)\leq n$, is $\rho$-structural.
 If $f \in \rho_0$, the identities
  \begin{equation}
    \label{eq:zerodim}
    f(\e_0,\e_2,\e_4,\e_6,\dots)=f(\e_0,\e_1, \ldots,\e_{n},\ldots)= f(\e_1,\e_3,\e_5,\e_7,\dots)
  \end{equation}
  specifying that $\dim(f)= 0$, are $\rho$-structural.
\end{example}

\begin{proposition} 
  \label{prop:str} 
  The class $\mathsf{Alg}_\rho^\top=\{\mathbf S^\top : \mathbf S\in \mathsf{Alg}_\rho\}$ is equal to $\Mod(\Str(\rho))$.
  \begin{proof}
    By Lemma \ref{lem:upterms}(i) every $\rho_\bullet$-algebra $\mathbf S^\top$ satisfies $\Str(\rho)$. 
    For the converse, let $\mathbf{A}$ be a $\rho_\bullet$-algebra such $\mathbf{A} \models \Str(\rho)$. 
    If $f \in \rho_\bullet$, there is $n \in \omega$ such that $f \in \rho_n$ and thus $\mathbf{A}$ satisfies \eqref{eq:findim} or \eqref{eq:zerodim}. 
    Let $\mathbf{S}$ be the $\rho$-algebra such that $f^{\mathbf{S}}(s_0, \ldots, s_{n-1})=f^{\mathbf{A}}(s_0, \ldots, s_{n-1}, s_0, s_0, \ldots)$. 
    Then $\mathbf{A}=\mathbf{S}^{\top}$. 
  \end{proof}
\end{proposition}

\begin{lemma}
  \label{lem_bullet}
  For a class $\mathsf{H}$ of $\rho$-algebras,
  $\Mod(\Th(\mathsf{H}^\top))=\Mod(\Str(\rho) \cup \Th(\mathsf{H})^\bullet)$. 
  \begin{proof}
    Let $\mathbf{A}$ be an $\rho_\bullet$-algebra that satisfies $\Str(\rho) \cup \Th(\mathsf{H})^\bullet$.
    Since $\mathbf{A} \models \Str(\rho)$, $\mathbf{A} = \mathbf{S}^{\top}$ for some $\rho$-algebra $\mathbf{S}$.
    We show that $\Th(\mathsf{H}^\top) \subseteq \Th(\mathbf{A})$. 
    Let $s = t \in \Th(\mathsf{H}^\top)$. 
    Then, by Lemma \ref{lem:upterms}(iii), $s^\circ = t^\circ \in \Th(\mathsf{H})$, so that, since $\mathbf{A}$ satisfies $\Th(\mathsf{H})^\bullet$, $\mathbf{A}$ satisfies $(s^\circ)^\bullet = (t^\circ)^\bullet$.
    By Lemma \ref{lem:upterms}(i), $\mathbf{A}$ satisfies $s = t$.
    Conversely, we show that $\Str(\rho) \cup \Th(\mathsf{H})^\bullet \subseteq \Th(\mathsf{H}^\top)$.
    That $\Str(\rho) \subseteq \Th(\mathsf{H}^\top)$ is clear by Proposition \ref{prop:str}.
    Let $u = v \in \Th(\mathsf{H})$. 
    Then, by Lemma \ref{lem:upterms}(iii) $u^\bullet = v^\bullet \in \Th(\mathsf{H}^\top)$ iff $(u^\bullet)^\circ = (v^\bullet)^\circ \in \Th(\mathsf{H})$. 
    The conclusion now follows from Lemma \ref{lem:upterms}(ii). 
  \end{proof}
\end{lemma}

\subsection{The main result} 
Now we use Theorem \ref{thm:bir1} and the last technical results to get 
an enhanced version of Birkhoff Theorem. 

To begin with, it is easy to see that the top operator commutes with $\mathbb{H}, \mathbb{S}, \mathbb{P}$. 

\begin{lemma}
  \label{lem_top}
  Let $\mathsf{H}$ be a class of $\rho$-algebras for some finitary $\rho$. 
  Let $\mathbf{S}, \mathbf{R}, \mathbf{S}_i (i \in I) \in \mathsf{H}$. 
  Then 
  \begin{itemize}
    \item $\prod_{i \in I} (\mathbf{S}_i^\top) = \left (\prod_{i \in I} \mathbf{S}_i \right ) ^{\top}$
    \item $\mathbf{R}$ is a subalgebra of $\mathbf{S}$ iff $\mathbf{R}^\top$ is a subalgebra of $\mathbf{S}^\top$; 
    \item $\alpha: \mathbf{R} \to \mathbf{S}$ is a homomorphism iff $\alpha: \mathbf{R}^{\top} \to \mathbf{S}^{\top}$ is a homomorphism. 
  \end{itemize}
  In particular $\mathbb{HSP}(\mathsf{H}^{\top})=(\mathbb{HSP}(\mathsf{H}))^{\top}$. 
\end{lemma}

\begin{theorem} 
  \label{thm_oldbirkh}
  The following are equivalent for a class $\mathsf{H}$ of $\rho$-algebras. 
  \begin{enumerate}
      \item $\mathsf{H}$ is a variety of $\rho$-algebras;
      \item $\mathsf{H}=\Mod(\Th(\mathsf{H}))$; 
      \item $\mathsf{H}^{\top}$ is a variety of $\rho_\bullet$-algebras;
      
      \item $\mathsf{H}^\top=(\mathsf{H}^\top)^{\uptriangle \downtriangle}$ and $(\mathsf{H}^\top)^\uptriangle$ is a variety of infinitary clone $\rho_{\bullet}$-algebras.
  \end{enumerate}
  \begin{proof}
      (1) $\IFF$ (3) By Lemma \ref{lem_top}. 

      \noindent
      (3) $\IMP$ (2) 
      Trivially, $ \mathsf{H} \subseteq \Mod(\Th(\mathsf{H}))$. 
      Conversely, assume that $\mathbf{S}$ satisfies $\Th(\mathsf{H})$. 
      The $\rho_{\bullet}$-algebra $\mathbf{S}^\top$ satisfies $\Str(\rho)$ by virtue of Proposition \ref{prop:str}. 
      Moreover, by Lemma \ref{lem:upterms}(iv), $\mathbf{S}^\top$ satisfies ${\Th(\mathsf{H})}^{\bullet}$. 
      Therefore $\mathbf{S}^\top \in \Mod(\Str(\rho) \cup {\Th(\mathsf{H})}^{\bullet})$, which, by Lemma \ref{lem_bullet}, is equal to $\Mod(\Th(\mathsf{H}^\top))$. 
      Then the conclusion follows from Theorem \ref{thm:bir1} by applying the hypothesis that $\mathsf{H}^{\top}$ is a variety of $\rho_\bullet$-algebras. 

      \noindent
      (2) $\IMP$ (1) 
      It is immediate to check that $\Mod(\Th(\mathsf{H}))=\mathsf{H}$ is a variety.

      \noindent
      (4) $\IFF$ (3) By Theorem \ref{thm:bir1}.
  \end{proof}
\end{theorem}

\section{Topological version of Birkhoff Theorem}
\label{sec:top}
The equivalence $(1) \IFF (2)$ of Theorem \ref{thm:bir1} can be interpreted to say, for a class $\mathsf{K}$ of $\tau$-algebras: $\mathbf{B}$ satisfies $\Th(\mathsf{K})$ iff $\mathbf{B} \in \mathbb{HSP}(\mathsf{K})$. 
When $\mathsf{K}=\{\mathbf{A}\}$, as the notion of satisfaction can rephrased in terms of homomorphisms of clones, the above reformulation gives the following.  
\begin{proposition}
  \label{cor_birk}
  Let $\mathbf{A}, \mathbf{B}$ be $\tau$-algebras. Then the following are equivalent: 
  \begin{enumerate}
    \item $\mathbf{B} \in \mathbb{HSP}(\mathbf{A})$; 
    \item the map $\epsilon: \mathbf{A}^{\uparrow} \to \mathbf{B}^{\uparrow}$, $\epsilon(t^{\mathbf{A}})=t^{\mathbf{B}}$, is a homomorphism of infinitary clone $\tau$-algebras. 
  \end{enumerate}
    \begin{proof}
     (1) $\IMP$ (2) It is immediate to check that if $\mathbf{B} \in \mathbb{P}(\mathbf{\mathbf{A}})$, $\mathbb{S}(\mathbf{\mathbf{A}})$, or $\mathbb{H}(\mathbf{\mathbf{A}})$, then $\epsilon$ is well-defined and it is a homomorphism. 
 
     \noindent
     (2) $\IMP$ (1) By virtue of Theorem \ref{thm:bir1} it is enough to prove that $\Th(\mathbf{A}) \subseteq \Th(\mathbf{B})$. 
     Now, since $\mathcal{N}_{\bar \tau}$ is initial we have the following factorisation 
     \begin{equation*}
       \begin{tikzcd}
           \mathcal{N}_{\bar \tau} \arrow[rr, "(-)^{\mathbf{B}}"] \arrow[dr, "(-)^{\mathbf{A}}", swap] & &  \mathbf{B}^{\uparrow}\\
           & \mathbf{A}^{\uparrow} \arrow[ur, "\epsilon", swap] & \\
       \end{tikzcd}
       \end{equation*}

       \vspace{-0.75cm}
   Observing that $\ker ((-)^\mathbf{A}) \subseteq \ker (\epsilon \circ (-)^\mathbf{A}) $, the above diagram implies that $\Th(\mathbf{A}) \subseteq \Th(\mathbf{B})$. 
   \end{proof}
\end{proposition}

Now we strenghten condition (2) in order to characterise $\mathbf{B} \in \mathbb{HSP}_{\fin}(\mathbf{A})$. 

\subsection{The topology of the pointwise convergence}
Let $A,B$ be sets. 
Consider, for $s\in A^\omega$ and $\phi \in O^{(\omega)}_A$ the sets
$\EuScript{U}_{s,\phi} :=\{\psi \in O^{(\omega)}_A : \psi(s)=\phi(s)\} \text{.}$
The family $\{\EuScript{U}_{s,\phi} : s \in A^\omega, \phi \in O^{(\omega)}_A\}$ forms a subbase for the \emph{topology of the pointwise convergence} on $O^{(\omega)}_A$. 
This is the product topology on $A^{A^\omega}$ regarding $A$ as a discrete space. 
We recall the notion of uniform continuous function between the spaces $O^{(\omega)}_A$ and $O^{(\omega)}_B$. 
Clearly a uniformly continuous function is always continuous, while the opposite is not generally the case. 

\begin{definition}
  \label{def_unif}
  Let $\epsilon: O^{(\omega)}_A \to O^{(\omega)}_B$ be a function. 
  We say that $\epsilon$ is 
 \emph{uniformly continuous} if for all $r \in B^\omega$ there are $s_0, \ldots, s_{n-1} \in A^\omega$ such that for every $\phi, \psi \in O^{(\omega)}_A$
  \begin{equation*}
    \phi(s_0) = \psi(s_0), \ldots,  \phi(s_{n-1}) = \psi(s_{n-1}) \IMP \epsilon(\phi)(r) = \epsilon(\psi)(r) \text{.}
  \end{equation*}
\end{definition}

Similarly, we can define the topology of the pointwise convergence on $O^{(n)}_A$, for each $n \in \omega$. 
The clone $O_A=\bigcup\{O_A^{(n)} : n \in \omega\}$ can be endowed with the disjoint union topology. 

\begin{definition}
  Let $\eta: O_A \to O_B$ be a function. 
  We say that $\eta$ is 
 \emph{uniformly continuous} if, for each $k \in \omega$, for all $r \in B^k$ there are $s_0, \ldots, s_{n-1} \in A^k$ such that for every $\phi, \psi \in O^{(k)}_A$
  \begin{equation*}
    \phi(s_0) = \psi(s_0), \ldots,  \phi(s_{n-1}) = \psi(s_{n-1}) \IMP \eta(\phi)(r) = \eta(\psi)(r) \text{.}
  \end{equation*}
\end{definition}

\subsection{Topological Birkhoff}
If $\mathbf{A}$ is a $\tau$-algebra for some homogeneous $\tau$, $\mathcal{O}^{(\omega)}_\mathbf{A}$ is a topological space with the pointwise convergence topology. 
Any functional infinitary clone $\tau$-algebra $\mathcal{F} \in \mathsf{FCA}(\mathbf{A})$ inherits this topology from $\mathcal{O}^{(\omega)}_\mathbf{A}$. 

\begin{theorem}
  \label{thm_topbirk}
  Let $\mathbf{A}, \mathbf{B}$ be $\tau$-algebras. Then the following are equivalent:
  \begin{enumerate}
    \item for every $r \in B^\omega$, $\mathbf{B}_{\bar r} \in \mathbb{HSP}_{\fin}(\mathbf{A})$; 
    \item the map $\epsilon: \mathbf{A}^{\uparrow} \to \mathbf{B}^{\uparrow}$, $\epsilon(t^{\mathbf{A}})=t^{\mathbf{B}}$, is a uniformly continuous
    homomorphism of infinitary clone $\tau$-algebras. 
  \end{enumerate}
  Moreover, if $\mathbf{A}$ and $\mathbf{B}$ are finite-dimensional, $(1)$ may be refined as: 
  \begin{enumerate}
    \item[($1'$)] for every $r \in B_{\fin}^\omega$, $\mathbf{B}_{\bar r} \in \mathbb{HSP}_{\fin}(\mathbf{A})$. 
  \end{enumerate}
  \begin{proof}
    (1) $\IMP$ (2) By Theorem \ref{thm_exp}, $\mathbb{HSP}(\mathbf{A})$ is closed under expansion, so that $\mathbf{B} \in \mathbb{HSP}(\mathbf{A})$. 
    Corollary \ref{cor_birk} implies that the map $\epsilon$ is a well-defined homomorphism. 
    We show that $\epsilon$ is uniformly continuous. 
    Let $r \in B^\omega$. 
    By Definition \ref{def_unif}, we need to find $s_0, \ldots, s_{n-1} \in A^\omega$ such that for every $t,u \in N_{\bar \tau}$ 
    \begin{equation*}
      t^{\mathbf{A}}(s_0) = u^{\mathbf{A}}(s_0), \ldots,  t^{\mathbf{A}}(s_{n-1}) = u^{\mathbf{A}}(s_{n-1}) \IMP t^{\mathbf{B}}(r) = u^{\mathbf{B}}(r) \text{.}
    \end{equation*}
    Let $\mathbf{D}$ be the subalgebra of $\mathbf{A}^n$ and $\alpha : \mathbf{D} \to \mathbf{B}_{\bar r}$ be the onto homomorphism that witnesses the assumption (1). 
    Since $\alpha$ is onto there are $s_0, \ldots, s_{n-1} \in A^\omega$ such that $\alpha^\omega(s_0, \ldots, s_{n-1})=r$. 
    If $t^{\mathbf{A}}(s_i) = u^{\mathbf{A}}(s_i)$ for each $i=0, \ldots, n-1$, since $\alpha$ is a homomorphism, 
    we have, as desired:  
    \begin{align*}
      t^{\mathbf{B}}(r) & = \alpha(t^{\mathbf{D}}(s_0, \ldots, s_{n-1})) \\
      & = \alpha(t^\mathbf{A}(s_0), \ldots, t^\mathbf{A}(s_{n-1}))  \\
      & = \alpha(u^\mathbf{A}(s_0), \ldots, u^\mathbf{A}(s_{n-1})) \\
      & = \alpha(u^{\mathbf{D}}(s_0, \ldots, s_{n-1})) \\
      & = u^{\mathbf{B}}(r)\text{.}
    \end{align*}

    \noindent
    (2) $\IMP$ (1) Let $r \in B^\omega$, and let $m \in (A^n)^\omega =(A^\omega)^n$ be witness that $\epsilon$ is uniformly continuous. 
    We show that there is an onto homomorphism $\alpha: (\mathbf{A}^n)_{\bar{m}} \to \mathbf{B}_{\bar r}$.
    Define 
    $\alpha(m_i)=r_i$ for all $i \in \omega$ and  
    $\alpha(t^{\mathbf{A}^{n}}(m)) = t^{\mathbf{B}}(r)$ for each $t \in N_{\bar \tau}$. 
    As $m \in (A^n)^\omega =(A^\omega)^n$, $m$ can be identified with $(s_0, \ldots, s_{n-1})$. 
    By uniform continuity $\alpha$ is well-defined, and by definition it is a surjective homomorphism. 

    \noindent
    ($1'$) $\IMP$ (2) By Lemma \ref{lem_top} $\mathbb{HSP}(\mathbf{A})$ is a class of finite-dimensional algebras.  
    By Proposition \ref{prop_exp}, $\mathbb{HSP}(\mathbf{A})$ is closed under finite expansion, so that $\mathbf{B} \in \mathbb{HSP}(\mathbf{A})$. 
    Then the map $\epsilon$ is a well-defined homomorphism. 
    The rest of the proof can be carried out as before.  
    
    \noindent
    (2) $\IMP$ ($1'$) Trivial. 
  \end{proof}
\end{theorem}

As a corollary, we get the following enhanced version of the topological Birkhoff theorem stated in 
\cite[Theorem 4.3]{GP18} (see also \cite{Sc17}). 

\begin{theorem}
  \label{thm:topbirk2}
  Let $\mathbf{R}$ and $\mathbf{S}$ be algebras over a finitary type $\rho$.
  The following are equivalent:
  \begin{enumerate}
    \item every finitely generated subalgebra of $\mathbf{S}$ belongs to $\mathbb{HSP}_{\fin}(\mathbf{R})$;
    \item the map $\eta: \Clo(\mathbf{R}) \to \Clo(\mathbf{S})$, $\eta(v^{\mathbf{R}})=v^{\mathbf{S}}$ for $v$ $k$-ary, is a uniformly continuous 
    homomorphism of clones; 
    \item the map $\epsilon: (\mathbf{R}^{\top})^{\uparrow} \to (\mathbf{S}^{\top})^{\uparrow}$, $\epsilon(t^{\mathbf{R}^{\top}})=t^{\mathbf{S}^{\top}}$, is a uniformly continuous
    homomorphism of infinitary clone $\rho_\bullet$-algebras.
  \end{enumerate}
  \begin{proof}
    (1) $\IFF$ (3) By Theorem \ref{thm_topbirk}. 

    \noindent
    (2) $\IFF$ (3) By surjectivity of the map $(-)^{\circ}$ and Lemma \ref{lem:upterms}. 
  \end{proof}
\end{theorem}

Finally, we note that with little adjustments in the proof, Theorem \ref{thm_topbirk} can lead to the following definition and result. 

\begin{definition}
  Let $\epsilon: O^{(\omega)}_A \to O^{(\omega)}_B$ be a function. 
  We say that $\epsilon$ is 
 \emph{globally uniformly continuous} if for all $r \in B^\omega$ there are $s_0, \ldots, s_n, \ldots \in A^\omega$ such that for every $\phi, \psi \in O^{(\omega)}_A$
  \begin{equation*}
    \phi(s_i) = \psi(s_i)   \text{ for all } i \in \omega \IMP \epsilon(\phi)(r) = \epsilon(\psi)(r) \text{.}
  \end{equation*}
\end{definition}

\begin{theorem}
  Let $\mathbf{A}, \mathbf{B}$ be $\tau$-algebras. Then the following are equivalent:
  \begin{enumerate}
    \item for every $r \in B^\omega$, $\mathbf{B}_{\bar r} \in \mathbb{HSP}_{\omega}(\mathbf{A})$; 
    \item the map $\epsilon: \mathbf{A}^{\uparrow} \to \mathbf{B}^{\uparrow}$, $\epsilon(t^{\mathbf{A}})=t^{\mathbf{B}}$, is a globally uniformly continuous
    homomorphism of infinitary clone $\tau$-algebras. 
  \end{enumerate}
  Moreover, if $\mathbf{A}$ and $\mathbf{B}$ are finite-dimensional, $(1)$ may be refined as: 
  \begin{enumerate}
    \item[($1'$)] for every $r \in B_{\fin}^\omega$, $\mathbf{B}_{\bar r} \in \mathbb{HSP}_{\omega}(\mathbf{A})$. 
  \end{enumerate}
\end{theorem}

\section{Conclusions and further work}
At the core of this work are the definition of infinitary clone algebras (Definition \ref{def:clo}) and their representation theorem (Proposition \ref{lem_neu}) adapted from \cite{neu70}. 
This definition revitalises Neumann approach, by adding to the type of $\aleph_0$-clones the operation symbols of a given type $\tau$ as nullary operations 
-- an idea already presented in \cite{BS22}.
From this apparatus stem the ``algebrisation of syntax'' discussed in Section \ref{sec:cats} and the two 
enriched versions of Birkhoff Theorem, both in the infinitary and the finitary case -- Theorems \ref{thm:bir1}, \ref{thm_oldbirkh}. 
Additionally, the significance of functional infinitary clone algebras is emphasized in the third, ``local'' version of Birkhoff Theorem, which addresses pseudovarieties of infinitary algebras (Theorem \ref{thm_topbirk}). 

We consider the techniques developed in this work as promising in order to revisit and enhance other 
classical results in universal algebra, 
and to provide new insights into open problems. 
In this direction, see for instance
the characterisation of the lattice of equational theories in terms of clone algebras presented in \cite[Theorem 5.7]{BS22}. 

Here are some lines for future research. 
\begin{itemize}
  \item To a given class of $\tau$-algebras $\mathsf{K}$, we can associate either $\mathsf{K}^{\uparrow}$, the pointwise image of $(-)^{\uparrow}$, or $\mathsf{K}^{\uptriangle}$. 
  In \ref{subs:central} we compared and distinguished their respective theories 
  thanks to the notion of central element (Lemma \ref{lemma_central}). 
  We say that a pure infinitary clone algebra is \emph{Boolean-like} if each of its elements is central. 
  The variety of infinitary Boolean-like algebras is axiomatized by the identities (N1)--(N3) and (C1)--(C3), and represents a generalisation of the variety of Boolean-like algebras of dimension $n$ \cite{SBLP20,BLPS18}. 
  This variety warrants further investigation.
 \item Consider the category whose objects are varieties of finitary algebras and whose morphisms are interpretations of varieties \cite[Section 10.8]{ALV3}.  
  This category, which is equivalent to the category of clones and clone homomorphisms, boasts a rich structure including a terminal and an initial object, $\mathsf{1}$ and $\mathsf{Set}$ respectively, 
  finite products, coproducts and tensor products, as defined in \cite[Definition 10.24]{ALV3}. 
  The ``decategorification'' of this category forms a bounded lattice. 
  Neumann \cite{neu74} noted that transitioning from a many-sorted to an infinitary one-sorted perspective in the treatment of clones -- transition which is advocated in the present article -- 
  may facilitate answering questions concerning this lattice and the interpretability of varieties. 
  Clone algebras, especially the infinitary clone algebras introduced here, provide a neat way to approach these problems.  
\end{itemize}

\noindent 
\emph{Data availability statement}. All data underlying the results are available as part of the article and no additional source data are required.

\bibliographystyle{abbrv}
\bibliography{birkhoff.bib}

\end{document}